\documentclass[10pt]{amsart}

\usepackage{latexsym,amssymb}

\usepackage{amsfonts,amsmath,amsthm,amsopn,amstext,amscd,xy,color,mathrsfs,etoolbox}

\input xy 
\xyoption{all}

\usepackage{url, enumitem}

\usepackage
[pdfauthor={Preston Wake and Carl Wang-Erickson},
 pdftitle={Deformation conditions for pseudorepresentations},
 bookmarks=false]
{hyperref}

\newtheorem{thm}{Theorem}[subsection] 
\newtheorem*{thm*}{Theorem}

\newtheorem{cor}[thm]{Corollary}
\newtheorem{prop}[thm]{Proposition}
\newtheorem{lem}[thm]{Lemma}

\theoremstyle{definition}
\newtheorem{defn}[thm]{Definition}
\newtheorem{defns}[thm]{Definitions}

\newtheorem{eg}[thm]{Example}

\theoremstyle{remark}
\newtheorem{rem}[thm]{Remark}

\newtheorem{warn}[thm]{Warning}

\makeatletter
\let\c@equation\c@thm
\makeatother
\numberwithin{equation}{subsection}

\newcommand\define{\newcommand}

\define\isoto{\xrightarrow{\sim}}
\define\onto{\twoheadrightarrow}

\DeclareMathOperator{\Spec}{Spec}

\define\cC{\mathcal{C}}

\newcommand{\ttmat}[4]{\left( \begin{array}{cc}
#1 & #2 \\
#3 & #4
\end{array}
\right)}

\newcommand{\Z}{\mathbb{Z}}
\newcommand{\Q}{\mathbb{Q}}

\newcommand{\F}{\mathbb{F}}

\newcommand{\m}{\mathfrak{m}}

\newcommand{\Hom}{\mathrm{Hom}}
\newcommand{\Gal}{\mathrm{Gal}}

\newcommand{\Ext}{\mathrm{Ext}}

\newcommand{\End}{\mathrm{End}}

\newcommand{\lb}{{[\![}}
\newcommand{\rb}{{]\!]}}

\newcommand{\red}{\mathrm{red}}

\define\cA{\mathcal{A}}

\define\GL{{\mathrm{GL}}}

\define{\Fitt}{\mathrm{Fitt}}
\define{\Ann}{\mathrm{Ann}}

\newcommand{\ra}{\rightarrow}

\newcommand{\lra}{\longrightarrow}
\newcommand{\lrisom}{\buildrel\sim\over\lra}
\newcommand{\risom}{\buildrel\sim\over\ra}
\newcommand{\rinj}{\hookrightarrow}

\newcommand{\rsurj}{\twoheadrightarrow}

\DeclareMathOperator{\Spf}{Spf}

\newcommand{\PsR}{{\mathrm{PsR}}}
\newcommand{\bF}{\mathbb{F}}
\newcommand{\bQ}{\mathbb{Q}}
\newcommand{\bZ}{\mathbb{Z}}
\newcommand{\bT}{\mathbb{T}}
\newcommand{\cE}{\mathcal{E}}

\newcommand{\Db}{{\bar D}}
\newcommand{\mDb}{\mathfrak m_{\bar D}}
\newcommand{\CRep}{{\mathcal R \mathrm{ep}}}
\newcommand{\Rep}{{\mathrm{Rep}}}
\newcommand{\Ad}{{\mathrm{GMA}}}

\newcommand{\cG}{{\mathcal{G}}}

\newcommand{\fl}{{\mathrm{flat}}}

\DeclareMathOperator{\Tr}{\mathrm{Tr}}
\DeclareMathOperator{\Sym}{\mathrm{Sym}}

\newcommand{\tr}{{\mathrm{tr}}}
\newcommand{\bro}{{\bar\rho}}

\newcommand{\sm}[4]{\ensuremath{\big(\begin{smallmatrix}#1 & #2 \\ #3 & #4\end{smallmatrix}\big)}}

\newcommand{\tfg}{{\mathrm{Tfg}}}
\newcommand{\CH}{{\mathcal{CH}}}

\newcommand{\Mt}{{\mathrm{Mod}^\mathrm{fin}_{\bZ_p[G]}}}
\newcommand{\PsDef}{{\mathrm{PsDef}}}

\setlist[enumerate]{label={(\arabic*)}}

\begin{document}

%

\makeatletter
\let\c@equation\c@thm
\makeatother
\numberwithin{equation}{subsection}

\title{Deformation conditions for pseudorepresentations}

\author{Preston Wake}
\address{Institute for Advanced Study\\
1 Einstein Drive \\
Princeton, NJ 08540}
\email{pwake@ias.edu}

\author{Carl Wang-Erickson}
\address{Department of Mathematics, Imperial College London \\
	London SW7 2AZ, UK}
\email{c.wang-erickson@imperial.ac.uk}

\setcounter{tocdepth}{1}
\begin{abstract}
Given a property of representations satisfying a basic stability condition, Ramakrishna developed a variant of Mazur's Galois deformation theory for representations with that property. We introduce an axiomatic definition of pseudorepresentations with such a property. Among other things, we show that pseudorepresentations with a property enjoy a good deformation theory, generalizing Ramakrishna's theory to pseudorepresentations.
\end{abstract}
\maketitle
 \tableofcontents


\section{Introduction}

The deformation theory of Galois representations has been used extensively to realize and study the proposed Langlands correspondence between Galois representations and automorphic representations. Depending on the setting in which the correspondence is being studied, one wants to set up the deformation theory so that it parameterizes exactly the Galois representations with certain properties. A collection of such properties is known as a ``deformation condition'' or ``deformation datum.'' Ramakrishna \cite{ramakrishna1993} proved that for any deformation condition $\cC$ satisfying a basic stability property, the deformation problem for $\cC$ is representable relative to the unrestricted deformation problem. 

In applications, and especially when working with residually reducible Galois representations, Galois pseudorepresentations are often more accessible than Galois representations. A Galois  pseudorepresentation is the data of a characteristic polynomial for each element of the Galois group, satisfying appropriate compatibility conditions. However, it has not been clear how to apply deformation conditions -- which apply to Galois modules -- to pseudorepresentations. 

The point of this paper is to solve this problem by axiomatically translating deformation conditions for representations into deformation conditions for pseudorepresentations. In particular, we construct universal deformation rings for pseudorepresentations satisfying such a condition. Also, it follows from our techniques that a deformation condition cuts out a closed locus in any family of Galois representations, generalizing the case of local families considered in \cite{ramakrishna1993}. We expect that the various constructions we make, especially the universal Cayley-Hamilton algebra defined by a condition, will find many applications. 

We were stimulated to develop this theory for application in the companion paper \cite{WWE3} (see also \cite{WWE5}). After this work was complete, the second-named author identified Galois cohomological data that controls the deformation theory of these pseudodeformation rings: see \cite{CarlAInf}, especially [Thm.\ 3.4.1, \textit{loc.\ cit.}]. 

\subsection{Background}
\label{subsec:background}

To give a more thorough explanation of our results, we define our scope. Fix a prime $p$ and a profinite group $G$. We assume that $G$ satisfies the finiteness condition $\Phi_p$ (see \S\ref{subsec:conv}). We are interested in studying the categories $\CRep_G$ of representations of $G$ and $\PsR_G$ of pseudorepresentations of $G$, and the related moduli spaces $\Rep_G$ and $\PsR_G$.

Let $\cC$ be a condition on finite-cardinality $\bZ_p[G]$-modules such that $\cC$ is closed under isomorphisms, subquotients, and finite direct sums (we will call such a condition \emph{stable}). Ramakrishna \cite{ramakrishna1993} showed that the formal deformation functor of representations with $\cC$ is representable. 

The goal of this paper is to sensibly extend the condition $\cC$ so that it applies to pseudorepresentations. One difficulty in doing this is that pseudorepresentations are defined as collections of functions; they need not come in any obvious way from a $\Z_p[G]$-module. 

The main tool we use to bridge the gap between $G$-pseudorepresentations and $G$-modules is the category $\CH_G$ of \emph{Cayley-Hamilton representations of $G$}, which was introduced by Chenevier \cite{chen2014}, building on work of Bella\"{i}che--Chenevier \cite{BC2009}, and further developed by second-named author \cite{WE2018}. The data of a Cayley-Hamilton representation of $G$ includes a $\Z_p[G]$-module, so it is possible to relate these objects to condition $\cC$. Moreover, $\CH_G$ is closely related to both $\CRep_G$ and $\PsR_G$; the key properties can be summarized as follows: 
\begin{itemize}[leftmargin=2em]
\item $\CH_G$ broadens the category of representations of $G$, in the sense that there is an inclusion $\CRep_G \rinj \CH_G$ and an induced-pseudorepresentation functor $\psi: \CH_G \ra \PsR_G$ extending the usual functor $\CRep_G \ra \PsR_G$.
\item $\psi$ is essentially surjective, which is not true for $\psi\vert_{\CRep_G}$. 
\item $\CH_G$ has a universal object, the ``universal Cayley-Hamilton representation'' of $G$, which we denote by $E_G$.
\item $E_G$ is finitely generated as a module over its base ring.
\end{itemize}
The first two conditions say that Cayley-Hamilton representations generalize the notation of representations to an extent sufficient to see all pseudorepresentations. The last two conditions mean that $\CH_G$ is still manageable in size; in particular, it reduces constructions to finite-cardinality $\Z_p[G]$-modules, so one can apply the condition $\cC$. We review the theory of Cayley-Hamilton representations in \S\ref{subsec:RepG}. 

\subsection{Results}
\label{subsec:results}
We use the above properties of $\CH_G$ to give the definition of a \emph{Cayley-Hamilton representation of $G$ with condition $\cC$ }(Definition \ref{defn:PforCH}). The Cayley-Hamilton representations of $G$ with condition $\cC$ form a full subcategory $\CH_G^\cC \subset \CH_G$. We say that a pseudorepresentation has condition $\cC$ if there is a Cayley-Hamilton representation that induces it. More precisely, we define the category of pseudorepresentations satisfying condition $\cC$, denoted by $\PsR_G^\cC$, as the essential image $\PsR_G^\cC \subset \PsR_G$ of $\psi$ restricted to $\CH_G^\cC \subset \CH_G$ (Definition \ref{defn:PsR-P}). 

We view these new definitions as being the heart of our paper. Our main results are meant to justify why we think this notion is a useful generalization of Ramakrishna's theory. We show that it has the same geometric and universality properties as Ramakrishna's, and that it specializes to his definition in cases where both definitions apply. Importantly, we also give tools for checking whether certain pseudorepresentations (including many of the ones that come up naturally) have this property. In some cases, we are also able to give tools for computing these objects in terms of natural Selmer-type groups.

Our main general results are as follows.
\begin{enumerate}[leftmargin=2em]
\item The category $\CH_G^\cC$ has a universal object $E_G^\cC$, which is a quotient of the universal object $E_G$ in $\CH_G$ (Theorem \ref{thm:univ_CH_P}).
\item We define the notion of a faithful Cayley-Hamilton $G$-module $M$ (it is a $G$-module with extra structures giving rise to a Cayley-Hamilton representation $E_M$ of $G$). We prove that $M$ satisfies $\cC$ if and only if $E_M$ satisfies $\cC$ (Theorem \ref{thm:ModToCH}). 
\item Property $\cC$ is a Zariski closed condition on both $\Rep_G$ and $\PsR_G$ (Corollary \ref{cor:P-closed}). In particular, for a fixed residual pseudorepresentation, there is a pseudodeformation-with-$\cC$ ring that is a quotient of the universal pseudodeformation ring (Theorem \ref{thm:RDP})
\end{enumerate}
Many of the pseudorepresentations that arise naturally in the study of modularity-lifting come from faithful Cayley-Hamilton $G$-modules (see Example \ref{eg:Tate module}), and result (2) is very useful for showing that such a pseudorepresentation has $\cC$. Result (3) is a generalization of a theorem of Ramakrishna \cite[Prop.~1.2]{ramakrishna1993} for residually irreducible representations, and, in that case, we prove that our construction agrees with his (Corollary \ref{cor:us=Ram}). 

In \S\S3-4, we restrict our attention to the subcategories of $\CH_G$, $\CRep_G$, and $\PsR_G$ that are \emph{residually multiplicity-free} (see Definition \ref{defn:RMF}). For the sake of this introduction, we abuse notation and use the same letters for the residually multiplicity-free versions. With this restriction, it is known that $E_G$ admits the extra structure of a \emph{generalized matrix algebra}, a theory developed by Bella\"{i}che--Chenevier (see \cite[\S1.5]{BC2009} and \cite[\S4.1]{WWE2}). Using this extra structure, we prove the following results in \S\ref{sec:GMAs} and \S\ref{sec:gma_ext}, respectively:
\begin{enumerate}[resume, leftmargin=2em]
\item The quotient $E_G^\cC$ of $E_G$ is uniquely characterized in terms of $\CRep_G^\cC$ (Theorem \ref{thm:GMA-P}). 
\item The structure of $E_G^\cC$ is related to extension groups with condition $\cC$ (Theorem \ref{thm:BCs and exts}). 
\end{enumerate}
We view this last result as especially important, as it allows one to compute with some of our constructions. In particular, it allows one to give a kind of upper bound on the size of the pseudodeformation ring  in terms of Selmer-type groups, and this is useful in proving modularity-lifting theorems (see Remark \ref{rem:use for R=T}). This method is used in our papers \cite{WWE3,WWE5}.

\begin{rem}
In the body of the paper, we mostly discuss variants of these results for the versions of the categories that have a fixed residual pseudorepresentation $\Db$. To translate into the form stated in this introduction, one can apply the results of \S \ref{subsec:RepG}, particularly Theorem \ref{thm:RDb}. 
\end{rem}

\subsection{Context}
We are motivated by the case that $G$ is a quotient of an absolute Galois group of a global field and $\cC = \{\cC_v\}$ is a condition on the restrictions of representations of $G$ to its decomposition groups $G_v$, especially those conditions $\cC_v$ that arise from $p$-adic Hodge theory when $v \mid p$. We discuss the examples we have in mind in \S\ref{sec:ex}.

The unavailability of a notion of pseudorepresentations satisfying said conditions has been an obstacle to generalizing the use of the deformation theory of Galois representations to the residually reducible case. In particular, in the papers \cite{WWE3,WWE5}, we require the particular condition $\cC_v$ that a $\bZ_p[G_v]$-module arises as the generic fiber of a finite flat group scheme over the ring of integers of a $p$-adic local field. (This is the same condition that especially motivated \cite{ramakrishna1993}.)

In our previous work \cite{WWE1, WWE2}, we considered the ordinary condition on 2-dimensional global Galois representations (see also \cite{CS2019} and \cite[\S7]{WE2018}). However, the ordinary condition is of a rather different flavor, as it does not apply readily to a finite cardinality $\bZ_p[G]$-module without extra structure.

\subsection{Acknowledgements} We thank K\c{e}stutis \v{C}esnavi\v{c}ius for pointing out an unusual discrepancy about the numbering in \cite{chen2014}. We thank James Newton for helpful comments and discussions about applications to modularity lifting and the paper \cite{tenauthors}, which led us to include Corollary \ref{cor:tenauthors}. We thank the anonymous referees for many helpful comments that have led to improvements in the exposition.

The intellectual debt owed to the works of Bella\"{i}che--Chenevier \cite{BC2009} and Chenevier \cite{chen2014} will be clear to the reader. We thank Jo\"el Bella\"{i}che for his interest in and encouragement about this work.

P.W.\ was supported by the National Science Foundation under the Mathematical Sciences Postdoctoral Research Fellowship No.\ 1606255 and and Grant No. DMS-1638352. C.W.E. was supported by the Engineering and Physical Sciences Research Council grant EP/L025485/1.

\subsection{Notation and conventions}
\label{subsec:conv}
Rings $R$ are commutative and Noetherian, unless otherwise noted. Algebras are associative but not necessarily commutative, and are usually finitely generated.  For an algebra $E$, the term ``$E$-module" is used to mean ``left $E$-module,'' unless otherwise stated, and we let $\mathrm{Mod}_E$ denote the category of left $E$-modules.

When $A$ is a local ring, to say that $B$ is commutative local $A$-algebra means that $B$ is a local ring equipped with a morphism of local rings $A \to B$.

We study integral $p$-adic representations and pseudorepresentations of a profinite group $G$, which is assumed to satisfy the $\Phi_p$ finiteness condition of Mazur \cite{mazur1989}, i.e.\ the maximal pro-$p$ quotient of every finite index subgroup of $G$ is topologically finitely generated. 

We will work with categories of topological rings $R$ discussed in \S\ref{subsec:RepG}, giving finitely generated $R$-modules a natural topology. Algebras $E$ over $R$ are also understood to have their natural topology, either being finitely generated as $R$-modules or of the form $R[G]$. Actions of $G$ and homomorphisms $G \ra E^\times$ are understood to be continuous. Pseudorepresentations $D : E \ra R$ (resp.\ $D : G \ra R$) are also understood to be continuous, which means that the characteristic polynomial coefficient functions $E \ra R$ (resp.\ $G \ra R$) defined in Definition \ref{defn:PsR}(8) are continuous. 

\subsubsection{A note on attribution} 
We wish to make it clear that the notions of pseudorepresentation, Cayley-Hamilton algebra, and generalized matrix algebra, which are used in an essential way in this paper, are due to Chenevier \cite{chen2014} and Bella\"iche--Chenevier \cite[Ch.\ 1]{BC2009}, building on works of others including Wiles \cite{wiles1988}, Procesi \cite{procesi1987}, Taylor \cite{taylor1991}, Nyssen \cite{nyssen1996}, and Rouquier \cite{rouquier1996}. 

Where possible, we have attempted to cite the original sources directly. In some cases, we also cite \cite{WE2018} for background on these objects because the category of profinite-topological Cayley-Hamilton representations of a profinite group was developed there, adapting the non-topological version of \cite[\S1.22]{chen2014}. In contrast, we note that Chenevier developed the profinite-topological notion of pseudorepresentation \cite[\S3]{chen2014}, which we rely on without making any additional contributions. 

\section{Deformation conditions for Cayley-Hamilton representations}
\label{sec:prop}

In this section, we develop the theory of Cayley-Hamilton representations with a deformation condition $\cC$. The first three subsections are mainly review. In \S \ref{subsec:compat}, we begin by recalling from \cite{chen2014} and \cite[\S2]{WE2018} the non-topological theory of pseudorepresentations and Cayley-Hamilton algebras, their compatible representations, and moduli spaces of these (see also \cite{WWE1} for further background and examples). In \S \ref{subsec:RepG}, we review the modifications of this theory for topological groups. In \S \ref{subsec:stability}, we discuss the stable conditions of Ramakrishna \cite{ramakrishna1993}, and show that there is a well-defined ``maximal quotient with $\cC$" of a finite cardinality $\Z_p[G]$-module.

The main new constructions are given in \S \ref{subsec:constr}. There, we show that the ``maximal quotient with $\cC$" of the universal Cayley-Hamilton algebra retains its algebra structure. We also show that, for an algebra quotient of a Cayley-Hamilton algebra, there is a maximal further quotient that is Cayley-Hamilton. Combining these constructions, we obtain a ``maximal Cayley-Hamilton quotient with $\cC$" of the universal Cayley-Hamilton algebra.

The main results and definitions of the section are stated and proven in \S \ref{subsec:CH-P}. We define the category of Cayley-Hamilton representations with $\cC$, and show that the object constructed in \S \ref{subsec:constr} is a universal object in this category. We define a pseudorepresentation with $\cC$ as a pseudorepresentation that comes from a Cayley-Hamilton representations with $\cC$, and prove that the corresponding deformation functor is pro-representable.

In \S \ref{subsec:ModCH-P}, we introduce the notion of Cayley-Hamilton $G$-module. Many pseudorepresentations that arise from modular forms (and more generally, automorphic forms) come from Cayley-Hamilton $G$-modules (see Example \ref{eg:Tate module}). We prove a criterion for when such pseudorepresentations have $\cC$. In \S \ref{subsec:formal moduli with C}, we apply this criterion to prove a generalization of Ramakrishna's theorem from formal families of representations to more general families.

\subsection{Pseudorepresentations, compatible representations, and Cayley-Hamilton algebras}
\label{subsec:compat}
Our starting point is following definition, due to Chenevier \cite{chen2014}.

\begin{defns}
\label{defn:PsR} 
Let $R$ be a ring, let $E$ be an $R$-algebra, and let $\cG$ a group. 
\smallskip
\begin{enumerate}[leftmargin=2em]
\item A \emph{pseudorepresentation}, denoted $D: E \ra R$ or $(E,D)$, is a multiplicative polynomial law\footnote{A polynomial law is a natural transformation $D: \underline{E \otimes_R (-)} \to \underline{(-)}$ of set-valued functors on commutative $R$-algebras, where $\underline{(-)}$ is the forgetful functor from $R$-algebras to sets. It is multiplicative if, for each commutative $R$-algebra $A$, the function $D_A: E \otimes_{R} A \to A$ satisfies $D_A(1)=1$, $D_A(xy)=D_A(x)D_A(y)$ for all $x, y \in E \otimes_R A$. It has degree $d$ if, for each commutative $R$-algebra $A$, we have $D_A(bx) = b^d D_A(x)$ for all $x \in E \otimes_R A$ and all $b \in A$.} of degree $d$ from $E$ to $R$, for some $d \ge 1$. We call $d$ the \emph{dimension} of $D$, and $R$ the \emph{scalar ring} of $(E,D)$. The data of $D$ consists of a function $E \otimes_R A \to A$ for each commutative $R$-algebra $A$, which we denote by $D_A$.

\item A \emph{pseudorepresentation of $\cG$ with coefficients in $R$}, denoted $D : \cG \ra R$, is a pseudorepresentation of $R[\cG]$.

\item If $D: E \to R$ is a pseudorepresentation, and $x \in E$, we define the \emph{characteristic polynomial} $\chi_D(x,t) \in R[t]$ by $\chi_D(x,t)=D_{R[t]}(t-x)$. It is monic of degree equal to the dimension $d$ of $D$. We define the \emph{trace} $\mathrm{Tr}_D(x)$ to be the additive inverse of the coefficient of $t^{d-1}$ in $\chi_D(x,t)$. 
\end{enumerate}
\end{defns}

\begin{rem}
This notion of pseudorepresentation is called a ``determinant'' by Chenevier. The notion of ``pseudorepresentation'' of a group was first considered by Wiles \cite[Lem.\ 2.2.3]{wiles1988} in the case $d=2$, and by Taylor \cite[\S 1]{taylor1991} in general. Taylor's definition was also considered by Rouquier \cite{rouquier1996}, who considered $R$-algebras (not just groups), and called the resulting objects ``pseudo-caract\`eres" (or ``pseudocharacters", in English). By \cite[Lem.\ 1.12]{chen2014}, $\mathrm{Tr}_D: E \to R$ is a ``pseudocharacter" in the sense of Taylor and Rouquier. By \cite[Prop.\ 1.29]{chen2014}, the map $D \mapsto \mathrm{Tr}_D$ is a bijection if $(2d)! \in R^\times$.
\end{rem}

\begin{eg}
\label{eg:det}
There is a $d$-dimensional pseudorepresentation $\det: M_d(R) \to R$ defined by letting $\det_A: M_d(A) \to A$ be the determinant for all commutative $R$-algebras $A$. More generally, if $E$ is an Azumaya $R$-algebra of degree $d$, e.g.\ $V$ is a projective $R$-module of rank $d$ and $E = \End_R(V)$, then there is a $d$-dimensional pseudorepresentation $\det: E \to R$ given by the reduced norm on $E$. 
\end{eg}

\begin{eg}
\label{eg:psi}
If $V$ is a projective $R$-module of rank $d$ and $\rho: R[\mathcal{G}] \to \End_R(V)$ is a representation of $\mathcal{G}$, then there is an associated pseudorepresentation $\psi(\rho):\mathcal{G}\to R$. It is defined, for a commutative $R$-algebra $A$ and $x \in A[\mathcal{G}]$, by 
\[
\psi(\rho)_A(x)=\mathrm{det}_A((\rho\otimes_R A)(x)),
\]
where $\det$ is as in the previous example. 
\end{eg}

Informally, $\psi(\rho)$ is the composition of $\rho$ and $\det$. In the following definition, we make this `composition' construction explicit, and introduce notation for some other standard ways to produce new pseudorepresentations from old ones.

\begin{defns}
\label{defn:NewPsFromOld} 
Let $f: R \to R'$ be a ring homomorphism, let $E$ and $E'$ be an $R$-algebra and an $R'$-algebra, respectively, and let $g: E\otimes_R R' \to E'$ be a morphism of $R'$-algebras. Let $\cG$ be a group. 
\smallskip
\begin{enumerate}[leftmargin=2em]
\item Let $D: E \to R$ be a pseudorepresentation. The \emph{base-change of $D$ by $f$}, denoted $f \circ D: E \otimes_R R' \to R'$, is the pseudorepresentation of $E\otimes_R R'$ defined by $(f \circ D )_A(x)=D_{A}(x)$, where $A$ is a commutative $R'$-algebra, and $x \in E\otimes_{R} A = (E \otimes_{R} R') \otimes_{R'} A$.

When $f$ is understood, we write $D \otimes_{R}R'$ instead of $f \circ D$.

\item Let $D': E' \to R'$ be a pseudorepresentation. The \emph{composition of $D'$ and $g$}, denoted $D' \circ g: E \otimes_R R' \to R'$, is the pseudorepresentation of $E\otimes_R R'$ defined by $(D' \circ g)_A(x)=D'_A((g\otimes_{R'} A)(x))$, where $A$ is a commutative $R'$-algebra and $x \in E\otimes_{R'} A$. 

\item A \emph{morphism of pseudorepresentations} $\rho:(E,D) \to (E',D')$ is the data of a pair $(f,g)$, such that $f \circ D= D' \circ g$. We define $\psi(\rho)= f \circ D= D' \circ g$.

A \emph{morphism of pseudorepresentations of $\cG$} is a morphism of pseudorepresentations $(R[\cG],D) \to (R'[\cG],D \otimes_R R')$ where the homomorphism $R'[\cG] \to R'[\cG]$ is the identity.

\item If $D: E \to R$ is a pseudorepresentation and $\rho: \cG \to E^\times$ is a homomorphism, there is an induced homomorphism $\tilde{\rho}: R[\cG] \to E$, which defines a morphism of pseudorepresentations $(R[\cG],D \circ \tilde{\rho}) \to (E,D)$. We abuse notation and denote this morphism  by $\rho$, and write $\psi(\rho):\cG \to R$ for the resulting pseudorepresentation of $\cG$.

\end{enumerate}
\end{defns}

It is not, in general, true that any pseudorepresentation $D:\mathcal{G} \to R$ is of the form $\psi(\rho)$ for some representation $\rho$ of $\mathcal{G}$ as in Example \ref{eg:psi}. It is, of course, true that $D=\psi(\rho)$ for some morphism $\rho$ of pseudorepresentations $(R[\mathcal{G}],D \circ \tilde{\rho}) \to (E,D)$ as in Definition \ref{defn:NewPsFromOld}(4) (we could take $E=R[\mathcal{G}]$ and $\tilde{\rho}$ to be the identity). However, it is a surprising fact that the analogous statement is true if we restrict $(E,D)$ to be in the more restrictive class of \emph{Cayley-Hamilton} pseudorepresentations {\cite[\S1.17]{chen2014}}: see Corollary \ref{cor: CH to PsR essentially surjective}.

\begin{defn}
	\label{defn:CH}
We call a pseudorepresentation $D : E \ra R$ \emph{Cayley-Hamilton} when $E$ is finitely generated as an $R$-algebra, and, for every commutative $R$-algebra $A$ and every $x \in E \otimes_R A$, the element $x$ satisfies the characteristic polynomial $\chi_D(x,t)\in A[t]$ associated to it by $D$. That is, $D$ is Cayley-Hamilton when $\chi_D(x,x) = 0$ for all $x \in E \otimes_R A$. When $D : E \ra R$ is Cayley-Hamilton, we call the pair $(E,D)$ a \emph{Cayley-Hamilton $R$-algebra}.
\end{defn}

An important property of Cayley-Hamilton algebras is the the following finiteness result.
\begin{prop} 
\label{prop:C-H alg is module-finite}
If $(E,D)$ is a Cayley-Hamilton $R$-algebra, then $E$ is finitely generated as an $R$-module. 
\end{prop}
\begin{proof}
This follows from \cite[Prop.\ 2.13]{WE2018}, which uses the theory of PI-algebras. Indeed, $E$ is finitely generated as an $R$-algebra (by definition of Cayley-Hamilton algebra) and $R$ is Noetherian by convention. 
\end{proof}

Of course, the pseudorepresentation $(M_d(R),\det)$ of Example \ref{eg:det} is a Cayley-Hamilton algebra (by the Cayley-Hamilton Theorem), and we think of Cayley-Hamilton representations as generalizations of this example. Just as algebra homomorphisms $E \to M_d(R)$ are given the special name \emph{representation}, we follow Chenevier and give certain morphisms of pseudorepresentations the special name \emph{Cayley-Hamilton representation}.

\begin{defns}
Let $\mathcal{G}$ be a group.
\smallskip
\begin{enumerate}[leftmargin=2em]	
\item A \emph{Cayley-Hamilton representation} of a pseudorepresentation $(E,D)$ is a morphism of pseudorepresentations $\rho:(E,D) \to (E',D')$ with $(E',D')$ a Cayley-Hamilton algebra.
\item A \emph{Cayley-Hamilton representation of $\cG$}, denoted $(E',\rho,D')$, is a morphism of pseudorepresentations $(R[\mathcal{G}],D) \to (E',D')$ such that the map $R[\mathcal{G}] \to E'$ comes from a homomorphism $\rho:\cG \to E'^\times$.
\item If $(E,D)$ is also Cayley-Hamilton, we also refer to a Cayley-Hamilton representation $(E,D) \to (E',D')$ as a \emph{morphism of Cayley-Hamilton algebras}.
\item A \emph{morphism of Cayley-Hamilton representations of $\cG$}, written $(E,\rho,D) \to (E',\rho',D')$, is a morphism of Cayley-Hamilton algebras $(E,D) \to (E',D')$ such that $\rho'=(E \to E') \circ \rho$.
\end{enumerate}
We let $\CH_\cG$ denote the category of Cayley-Hamilton representations of $\cG$, with morphisms as defined in (3).
\end{defns}

\begin{eg}
\label{eg:compat_is_CH}
Let $(E,D)$ be a pseudorepresentation over $R$, and let $V$ be a projective $A$-module of rank $d$ for a commutative $R$-algebra $A$. Let $(\End_A(V), \det)$ be the pseudorepresentation defined in Example \ref{eg:det}, which is a Cayley-Hamilton algebra. A Cayley-Hamilton representation $(E,D) \to (\End_A(V),\det)$ is the data of a representation $\rho: E \to \End_A(V)$ of $E$ that is compatible with the pseudorepresentation, i.e.~such that $D=\det \circ \rho$. 
\end{eg}

This example motivates the following (perhaps more familiar) definition.
\begin{defn}
\label{defn:compat_functor}
Let $D : E \ra R$ be a a $d$-dimensional pseudorepresentation. Let $A$ be a commutative $R$-algebra. 

\begin{enumerate}[leftmargin=2em]
\item A \emph{framed compatible representation of $(E,D)$ over $A$} is an $R$-algebra homomorphism $\rho_A: E \ra M_d(A)$ such that $D \otimes_R A = \det \circ \rho_A$. We denote by $\Rep^\square_{E,D}$ the set-valued functor which assigns to a commutative $R$-algebra $A$ the set of all framed compatible representations of $(E,D)$ over $A$. 

\item A \emph{compatible representation of $(E,D)$ over $A$} is a pair $(V_A, \rho_A)$, where $V_A$ is a projective rank $d$ $A$-module and $\rho_A : E \ra \End_A(V)$ is an $R$-algebra homomorphism such that $D \otimes_R A = \det \circ \rho_A$. We denote by $\Rep_{E,D}$ the $\Spec R$-groupoid (i.e.\ stack) which assigns to a commutative $R$-algebra $A$ a category whose objects are the compatible representations of $(E,D)$ over $A$ and whose morphisms are isomorphisms of such data. 
\end{enumerate}
\end{defn}

\begin{rem}
The functor $\Rep^\square_{E,D}$ and the groupoid $\Rep_{E,D}$ are representable over $\Spec R$; see e.g.\ Theorem \ref{thm:ED-compat}. In particular, the natural adjoint action of $\GL_d$ on the affine $\Spec R$-scheme $\Rep^\square_{E,D}$ provides a smooth presentation for $\Rep_{E,D}$ as a $\Spec R$-algebraic stack. 
\end{rem}

Note that, by Example \ref{eg:compat_is_CH}, compatible representations are a special kind of Cayley-Hamilton representation. We view Cayley-Hamilton representations as a natural generalization of compatible representations.

\subsection{Representations of a profinite group}
\label{subsec:RepG}
For the remainder of the section, we fix a profinite group $G$ that satisfies the $\Phi_p$ condition. Recall the conventions from \S \ref{subsec:conv} regarding continuity.

Let $\mathrm{Adm}_{\Z_p}$ denote the category of admissible topological (not necessarily Noetherian) rings  where $p$ is topologically nilpotent (see \cite[Def.~0.7.1.2, pg.~60]{ega1} for the definition of admissible). We let $\tfg_{\Z_p} \subset\mathrm{Adm}_{\Z_p}$ be the full subcategory of topologically finitely generated objects (i.e.~those $A \in \mathrm{Adm}_{\Z_p}$ for which there exists a (non-topological) homomorphism $\Z_p[x_1,\dots,x_n] \to A$ with dense image), which are Noetherian rings. For $R \in \tfg_{\Z_p}$, we define $\tfg_R$ as the slice category of $\tfg_{\Z_p}$ over $R$. We often use the following fact.

\begin{lem}
\label{lem:tfg algebra} 
Let $R \in \tfg_{\Z_p}$ and let $A \in \tfg_R$. For any $x \in A$ non-zero, there exists an commutative local $R$-algebra of finite cardinality and an $R$-algebra homomorphism $f:A \to B$ such that $f(x) \ne 0$. In particular, a surjection $g:A \to A'$ in $\tfg_R$ is determined by the natural transformation $\Hom_{R-\mathrm{alg}}(A',-) \to \Hom_{R-\mathrm{alg}}(A,-)$ of functors on commutative local $R$-algebras of finite cardinality. 
\end{lem}
\begin{proof}
We leave this as an exercise. The main point is that, for any ideal of definition $I$ for $A$, the ring $A/I^n$ is a finitely generated $\Z$-algebra, and consequently a Noetherian Jacobson ring. 
\end{proof}

\begin{defn}
	For $A \in \tfg_{\Z_p}$, a \emph{representation of $G$ with coefficients in $A$} is a finitely generated projective $A$-module $V_A$ of constant rank with an $A$-linear $G$-action $\rho_A$.	We write $\Rep_G$ for the category of representations, fibered in groupoids over $\tfg_{\Z_p}$ via the forgetful functor $(V_A, \rho_A) \mapsto A$. We write $\Rep_G^\square$ for the category defined just as $\Rep_G$, but with the additional data of an $A$-basis for $V_A$. 
	
	Write $\PsR_G$ for the category of pseudorepresentations of $G$, fibered in groupoids over $\tfg_{\Z_p}$. 
	
	We write $\psi$ for the functor $\psi: \Rep_G \ra \PsR_G$ over $\tfg_{\Z_p}$ that sends a representation $\rho_A$ to its induced pseudorepresentation $\psi(\rho_A)$.
	
	We decompose each of these categories by dimension $d \geq 1$, writing $\Rep_G^d \subset \Rep_G$, etc.
\end{defn}

To understand $\PsR_G^d$, we introduce deformations. Given a finite field $\F/\F_p$, we let $\hat\cC_{W(\F)} \subset \tfg_{\Z_p}$ be the category of complete Noetherian commutative local $W(\bF)$-algebras $(A, \m_A)$ with residue field $\bF$. 

\begin{defn}
Let  $\Db: G \ra \bF$ be a pseudorepresentation. Its deformation functor $\PsDef_{\Db} : \hat\cC_{W(\F)} \ra \mathrm{Sets}$ is 
\begin{equation}
\label{eq:PsR_functor}
A \mapsto \{D : A[G] \ra A \text{ such that } D \otimes_A \bF \simeq \Db\}
\end{equation}
and elements of $\PsDef_{\Db}(A)$ are called \emph{deformations of $\Db$} or \emph{pseudodeformations}. 
\end{defn}

Remarkably, when one varies $\Db$ over a certain set of finite field-valued pseudorepresentations known as \emph{residual pseudorepresentations} (see \cite[Def.\ 3.4]{WE2018} for the definition), $\PsDef_\Db$ captures all of $\PsR_G^d$, in the following sense. 

\begin{thm}[Chenevier]
\label{thm:RDb}
Assume $G$ satisfies $\Phi_p$. 
\begin{enumerate}[leftmargin=2em]
	\item Given a finite field $\bF$ and a pseudorepresentation $\Db : G \ra \bF$, the functor $\PsDef_\Db$ is represented by an object  $(R_\Db, \mDb) \in \hat{\cC}_{W(\F)}$. 
	\item There is an isomorphism 
	\begin{equation}
	\label{eq:PsR-decomp}
	\PsR_G^d \cong \coprod_\Db \Spf R_\Db,
	\end{equation}
	where $\Db$ varies over $d$-dimensional residual pseudorepresentations. 
\end{enumerate}
\end{thm}

\begin{proof}
See \cite[Prop.\ 3.3, Prop.\ 3.7 and Cor.\ 3.14]{chen2014}.
\end{proof}

Let $\Rep_\Db$ (resp.\ $\Rep_\Db^\square$) denote the fiber in $\Rep_G$ (resp.\ $\Rep_G^\square$) of $\psi$ over $\Spf R_\Db$, where $\Db$ is a residual pseudorepresentation, so that we have
\begin{equation}
\label{eq:Rep-decomp}
\Rep_G^d \cong \coprod_\Db \Rep_\Db, \qquad \Rep_G^{\square,d} \cong \coprod_\Db \Rep_\Db^\square
\end{equation}
where $\Db$ varies over $d$-dimensional residual pseudorepresentations.

Now, and for the rest of this section, we fix a finite field $\bF/\bF_p$ and a pseudorepresentation $\Db:G \to \F$. In light of the decompositions \eqref{eq:PsR-decomp} and \eqref{eq:Rep-decomp}, we lose no scope in our study of $p$-adic families by fixing this choice.

\begin{defn}
	Let $A \in \tfg_{\Z_p}$. We say that a pseudorepresentation $D : G \ra A$ \emph{has residual pseudorepresentation $\Db$} when $\Spf A \ra \PsR_G^d$ is concentrated over $\Spf R_\Db$. We write $D^u_\Db : G \ra R_\Db$ for the universal pseudodeformation of $\Db$. 
	 
	A Cayley-Hamilton representation $(E, \rho: G \to E^\times, D: E \ra A)$ of $G$ over $A \in \tfg_{\bZ_p}$ \emph{has residual pseudorepresentation $\Db$} if its induced pseudorepresentation $\psi(\rho) : G \ra A$ has residual pseudorepresentation $\Db$. In particular, a representation $(V_A, \rho_A) \in \Rep_G^d(A)$ has residual pseudorepresentation $\Db$ when $\psi(\rho_A)$ does. 
	
	We let $\CH_{G,\Db}$ denote the full subcategory of $\CH_G$ (introduced in \S\ref{subsec:background}) whose objects have residual pseudorepresentation $\Db$. The natural transformation 
	\begin{equation}
	\label{eq: Rep to CH-Rep}
	\Rep_\Db \ra \CH_{G,\Db}
	\end{equation}
	arises from the Cayley-Hamilton representation structure on a representation of Example \ref{eg:psi}. 
	This natural transformation commutes with the induced pseudorepresentation functor $\psi$. 
\end{defn}

We observe that $\Rep_\Db$ parameterizes representations of $G$ with residual pseudorepresentation $\Db$. From \cite[Thm.\ 3.8]{WE2018} we know that $\Rep_\Db$ and $\Rep^\square_\Db$ admit natural algebraic models over $\Spec R_\Db$. By this we mean that there exists a finite type affine scheme (resp.\ algebraic stack) over $\Spec R_\Db$ whose $\m_\Db$-adic completion is $\Rep^\square_\Db$ (resp.\ $\Rep_\Db$). This implies that $\Rep^\square_\Db$ is an affine formal scheme and $\Rep_\Db$ is a formal algebraic stack, both formally of finite type over $\Spf R_\Db$. 

This algebraic model arises from a canonical universal Cayley-Hamilton representation of $G$ with residual pseudorepresentation $\Db$, which we now define. 

\begin{thm}
	\label{thm:univ_CH}
	The category $\CH_{G,\Db}$ has a universal object $(E_\Db, \rho^u: G \ra E_\Db^\times, D^u_{E_\Db} : E_\Db \ra R_\Db)$. In particular, $E_\Db$ is a finitely generated $R_\Db$-module. The map $\rho^u : R_\Db[G] \ra E_\Db$ is surjective and $D^u_{E_\Db} : E_\Db \ra R_\Db$ is a factorization of the universal pseudodeformation $D^u_\Db : G \ra R_\Db$ through $E_\Db$. 
\end{thm}

\begin{proof}
	See \cite[Prop.\ 3.6]{WE2018}. 
\end{proof}

Another way of stating the final sentence of Theorem \ref{thm:univ_CH} is that the induced pseudorepresentation $\psi(\rho^u) := D^u_{E_\Db} \circ \rho^u : G \ra R_\Db$ of the universal Cayley-Hamilton representation $\rho^u$ (with residual pseudorepresentation $\Db$) is equal to the universal deformation $D^u_\Db : G \ra R_\Db$ of $\Db$. 

The following theorem is proved using \eqref{eq: Rep to CH-Rep} and the universal property of $\rho^u$. 

\begin{thm}[{\cite[Thm.\ 3.7]{WE2018}}]
\label{thm:ED-compat}
There is an isomorphism of topologically finite type formal algebraic stacks (resp.\ formal schemes) on $\tfg_{\bZ_p}$, 
\[
\Rep_\Db \lrisom \Rep_{E_\Db, D^u_{E_\Db}}, \qquad \Rep^\square_\Db \lrisom \Rep^\square_{E_\Db, D^u_{E_\Db}}.
\]
\end{thm}

\begin{cor}
\label{cor: CH to PsR essentially surjective}
Let $A \in \tfg_{\Z_p}$. Let $D : G \ra A$ be a pseudorepresentation. Then there exists a Cayley-Hamilton algebra over $A$ that induces it. 
\end{cor}
\begin{proof}
By Theorem \ref{thm:RDb}, we reduce to the case that $D$ has a fixed residual pseudorepresentation and deduce that there exists a homomorphism $R_\Db \ra A$. Then the Cayley-Hamilton representation 
\[
(E_\Db \otimes_{R_\Db} A, \rho^u \otimes_{R_\Db} A : G \ra E_\Db \otimes_{R_\Db} A, D^u_{E_\Db} \otimes_{R_\Db} A)
\]
of $G$ over $A$ induces $D$. 
\end{proof}

\subsection{Stability}
\label{subsec:stability}
 
Let $\Mt \subset \mathrm{Mod}_{\Z_p[G]}$ be the full subcategory whose objects have finite cardinality.

\begin{defn}
\label{defn:stable} A \emph{condition} on $\Mt$ is a full subcategory $\cC \subset \Mt$. We will say an object of $\Mt$ \emph{satisfies condition $\cC$} or \emph{has $\cC$} if the object is in $\cC$.

A condition $\cC$ is \emph{stable} if it is preserved under isomorphisms, subquotients, and finite direct sums in $\Mt$. In other words, $\cC$ is stable if
\begin{enumerate}
\item for every object $A$ in $\cC$ and all isomorphisms $f:A \to B$ in $\Mt$, the object $B$ is also in $\cC$, and
\item for every object $A$ in $\cC$ and all morphisms $f:A \to B$ and $g:C \to A$ in $\Mt$, the kernel of $f$ and cokernel of $g$ are in $\cC$, and
\item for every finite collection of objects $A_1, \dots, A_n$ of $\cC$, the direct sum $A_1 \oplus \dots \oplus A_n$ in $\Mt$ is an object of $\cC$.
\end{enumerate}
\end{defn}

\begin{eg}
\label{eg:factors through is stable}
Let $H \subset G$ be a normal subgroup. Let $\cC \subset \Mt$ be the full subcategory of objects where the $G$ action factors through the quotient $G \to G/H$. Then $\cC$ is stable (as $\cC \cong \mathrm{Mod}_{\Z_p[G/H]}^\mathrm{fin}$).
\end{eg}

\begin{eg}
\label{eg:fiber product of stable is stable}
Let $H_1,\dots,H_n \subset G$ be subgroups, and, for $i=1,\dots,n$, let $\cC_i \subset \mathrm{Mod}_{\Z_p[H_i]}^{\mathrm{fin}}$ be a condition. Then there is a condition $\cC \subset \Mt$ defined by the Cartesian square
\[
\xymatrix{
\cC \ar[r] \ar[d] & \prod_{i=1}^n \cC_i \ar[d] \\
\Mt \ar[r] & \prod_{i=1}^n \mathrm{Mod}_{\Z_p[H_i]}^{\mathrm{fin}}.
}
\]
In other words, an object $M \in \Mt$ has $\cC$ if and only if the restriction $M|_{H_i} \in \mathrm{Mod}_{\Z_p[H_i]}^{\mathrm{fin}}$ has $\cC_i$ for all $i=1,\dots,r$. If all $\cC_i$ are stable, then this $\cC$ is stable.
\end{eg}

For examples of conditions $\cC$ that are of use in arithmetic, see \S\ref{sec:ex}. For the rest of this section, we fix a stable condition $\cC$.

\begin{thm}[Ramakrishna]
\label{thm:Ram}
	Let $A$ be a complete commutative Noetherian local $\bZ_p$-algebra and let $V_A$ be a finitely generated free $A$-module with an $A$-linear $G$-action. Then there exists a maximal quotient $A \rsurj A^\cC$ such that for an commutative local $A$-algebra $B$ of finite cardinality, the $\Z_p[G]$-module $V_A \otimes_A B$ satisfies $\cC$ if and only if $A \ra B$ factors through $A^\cC$.
\end{thm}
\begin{proof}
This follows immediately from \cite[Thm.\ 1.1]{ramakrishna1993}.
\end{proof}

\begin{lem}
\label{lem:Cstart}
The inclusion functor $\cC \to \Mt$ has a left adjoint $(-)^\cC:\Mt \to \cC$. For $V \in \Mt$, the functor $V \mapsto V^\cC$ has the following additional properties:
\begin{enumerate}
\item there is a quotient map $f_V: V \rsurj V^\cC$ of $\Z_p[G]$-modules.
\item for $W \in \cC$, the adjunction isomorphism
\[
\Hom_\cC(V^\cC,W) \cong \Hom_\Mt(V,W)
\]
is given by $g \mapsto g\circ f_V$.
\item if $A$ is a commutative $\bZ_p$-algebra and $V \in \mathrm{Mod}_{A[G]}$ as well, then $V^\cC \in \mathrm{Mod}_{A[G]}$. 
\end{enumerate}
\end{lem}

\begin{proof}
Let $V \in \Mt$. First we construct the quotient map $f_V: V \to V^\cC$. Since $V$ has finite cardinality, there are a finite number of quotients of $V$. Let $\{V_1, \dots, V_n\}$ be the (possibly empty) set all of quotients of $V$ that have property $\cC$. Define $f_V: V \to V^\cC$ to be the quotient by the kernel of $V \to \oplus_{i=1}^n V_i$. Then $V^\cC$ is isomorphic to a submodule of $\oplus_{i=1}^n V_i$. Since $\cC$ is closed under isomorphisms, subobjects, and finite direct sums, we see that $V^\cC$ satisfies $\cC$. By definition, any of the quotients $V \onto V_i$ factor through $f_V$, and this factoring is unique since $f_V$ is an epimorphism.

We now let $W \in \cC$ and show that the map given in (2)\ is an isomorphism. Let $W \in \cC$ and let $f:V \to W$ be a homomorphism of $\Z_p[G]$-modules. Let $W'=V/\ker(f)$, let $f':V \to W'$ be the quotient map, and let $f'':W' \to W$ be the injection induced by $f$. Since $W'$ is isomorphic to a submodule of $W$, we have $W' \in \cC$. Then $f':V \to W'$ must be one of the quotient maps $V \to V_i$, so $f'$ factors uniquely through $f_V$, i.e.~there is a unique morphism $g': V^\cC \to W'$ such that $f'=g' \circ f_V$. Now let $g:V^\cC \to W$ be $g=f''  \circ g'$; a simple computation shows that this assignment $f \mapsto g$ is inverse to the map given in (2).

Now we show that $V \mapsto V^\cC$ is a functor. Let $h: V \to V'$ be a $\Z_p[G]$-module homomorphism. Then $f_{V'}\circ h \in \Hom_\Mt(V,(V')^\cC)$, so by (2)\ there is a unique map $h^\cC: V^\cC \to (V')^\cC$ such that $f_{V'}\circ h=h^\cC \circ f_V$. This shows that $V \mapsto V^\cC$ is a functor, and, together with (2),\ this implies that it is left adjoint to the inclusion.

Finally, (3)\ follows from the functoriality. Indeed, let $a \mapsto m_a$ denote the structure map $A \to \End(V)$. Then there is a unique map $m^\cC_{\,\cdot}:A \to \End(V^\cC)$ such that $m^\cC_a \circ f_V = f_V \circ m_a$ for all $a \in A$. The fact that $m^\cC_{\,\cdot}$ is a ring homomorphism follows from the fact that $f_V$ is an epimorphism; for example, if $a,a'\in A$, we have 
\[
m^\cC_{aa'} \circ f_V=f_V \circ m_{aa'} = f_V \circ m_a \circ m_{a'}
\]
and
\[
m^\cC_{a} \circ m^\cC_{a'} \circ f_V = m^\cC_{a} \circ f_V \circ m_{a'} = f_V \circ  m_{a}  \circ m_{a'},
\]
from which we conclude $m^\cC_{aa'}=m^\cC_{a} \circ m^\cC_{a'}$.
\end{proof}

\subsection{Constructions}
\label{subsec:constr}
In what follows, we will treat finite cardinality left $E_\Db$-modules as objects of $\Mt$ via the map $\Z_p[G] \ra E_{\Db}$. By Proposition \ref{prop:C-H alg is module-finite}, $E_\Db$ is finite as a $R_\Db$-module. In particular, for any $i \ge 1$, the $\Z_p[G]$-module $E_\Db/\mDb^i E_\Db$ has finite cardinality, and thus is an object of $\Mt$.

\begin{lem}
\label{lem:max-cC-quot-E_Db}
For any $i\ge 1$, there is a unique $E_\Db$-module quotient $E_\Db/\mDb^i E_\Db \onto E^\cC_\Db(i)$ such that $E^\cC_\Db(i)$ has $\cC$ and such that any $E_\Db$-module quotient $E_\Db/\mDb^i E_\Db \onto W$ where $W$ has $\cC$ factors uniquely through $E_\Db/\mDb^i E_\Db \onto E^\cC_\Db(i)$.

\end{lem}
\begin{proof}
By Theorem \ref{thm:univ_CH}, $R_\Db[G] \ra E_\Db$ is surjective, so the lattice of $E_\Db$-quotients and $R_\Db[G]$-quotients of $E_\Db/\mDb^i E_\Db$ are identical. By Lemma \ref{lem:Cstart}, we can take $E^\cC_\Db(i) = (E_\Db/\mDb^i E_\Db)^\cC$.
\end{proof}

\begin{lem}
\label{lem:E(i)-compat}
For any $i \geq 1$, there is a canonical isomorphism $E^\cC_\Db(i+1) \otimes_{R_\Db} R_\Db/\mDb^i \risom E^\cC_\Db(i)$. 
\end{lem}
\begin{proof}
Let $E'=E^\cC_\Db(i+1) \otimes_{R_\Db} R_\Db/\mDb^i$. We show that $E'$ satisfies the universal property of Lemma \ref{lem:max-cC-quot-E_Db}. The composite 
\[
E_\Db/\mDb^{i+i} E_\Db \onto E^\cC_\Db(i+1) \onto E'
\]
factors through $E_\Db/\mDb^i E_\Db$, so $E'$ is a quotient of $E_\Db/\mDb^i E_\Db$. Since $E^\cC_\Db(i+1)$ has $\cC$ and $E'$ is a quotient of it, we see that $E'$ has $\cC$.

Now suppose that $E_\Db/\mDb^i E_\Db \onto W$ where $W$ has $\cC$. By the universal property of $E^\cC_{\Db}(i+1)$, the composite
\[
E_\Db/\mDb^{i+1} E_\Db \onto E_\Db/\mDb^i E_\Db \onto W
\]
factors uniquely through a map $E^\cC_{\Db}(i+1) \onto W$. Since $W$ is a $R/\mDb^i$-module, this factors uniquely through $E' \onto W$.
\end{proof}

\begin{lem}
\label{lem:EP}
For any $i \ge 1$, the module quotient $E_\Db/\mDb^i E_\Db \onto E^\cC_\Db(i)$ has the following properties. 
\begin{enumerate}[leftmargin=2em]
\item Let $N$ be a left $E_\Db/\mDb^i E_\Db$-module that has finite cardinality. Then $N$ satisfies condition $\cC$ as a $\bZ_p[G]$-module if and only if every map of left $E_\Db/\mDb^i E_\Db$-modules $E_\Db/\mDb^i E_\Db \ra N$ factors through $E^\cC_\Db(i)$. 
\item There is a natural right action of $E_\Db/\mDb^i E_\Db$ on $E^\cC_\Db(i)$, making $E^\cC_\Db(i)$ a quotient $R_\Db$-algebra of $E_\Db$.
\item An $E_\Db/\mDb^i E_\Db$-module $N$ that has finite cardinality satisfies condition $\cC$ if and only if its $E_\Db$-action factors through $E^\cC_\Db(i)$. 
\end{enumerate}
\end{lem}
\begin{proof}
(1)\ If $N$ satisfies $\cC$, then so does the image of any $E_\Db/\mDb^i E_\Db \ra N$, so this arrow factors through $E^\cC_\Db(i)$ by Lemma \ref{lem:max-cC-quot-E_Db}. Conversely, if every such arrow factors through $E^\cC_\Db(i)$, then because of the finiteness assumption on $N$ there exists some $m \in \bZ_{\geq 1}$ and a surjective map $(E^\cC_\Db(i))^{\oplus m} \rsurj N$. Consequently $N$ satisfies $\cC$. 

(2)\ Choose $z \in E_\Db$ and consider the composite morphism of left $E_\Db$-modules (they are also in $\Mt$) 
\[
E_\Db/\mDb^i E_\Db \buildrel{(\ ) \cdot z}\over\lra E_\Db/\mDb^i E_\Db \lra E^\cC_\Db(i)
\]
where the leftmost arrow is right multiplication by $z$. The composite must factor through $E^\cC_\Db(i)$ by (1). The resulting map $E^\cC_\Db(i) \ra E^\cC_\Db(i)$ gives the desired right action of $z$ on $E^\cC_\Db(i)$ and shows that $E^\cC_\Db(i)$ is an $R_\Db$-algebra. 

(3)\ This follows directly from (1) and (2) in light of the following general fact: for an algebra $E$ and a left $E$-module $M$, the $E$-action on $M$ factors through a quotient algebra $E \rsurj Q$ if and only if every morphism of left $E$-modules $E \ra M$ factors through $Q$. This follows from the fact that any such $E \ra M$ is of the form $x \mapsto x \cdot m$ for some $m \in M$.
\end{proof} 

By Lemmas \ref{lem:E(i)-compat} and \ref{lem:EP}(2), we have an inverse system $\{E_\Db^\cC(i)\}$ of $R_\Db$-algebra quotients of $E_\Db$. In particular, we have the $R_\Db$-algebra quotient 
\begin{equation}
\label{eq:KP}
E_\Db \rsurj \varprojlim_i E^\cC_\Db(i).
\end{equation}

This quotient almost gives the algebra we are looking for. However, we need to produce, from this algebra quotient, a Cayley-Hamilton algebra. The following lemma gives the general procedure for doing this. A special case of this lemma was employed in \cite{WWE1,WWE2}.

\begin{lem}
\label{lem:max_CH_quot}
Let $(E, D : E \ra R)$ be a Cayley-Hamilton $R$-algebra. Let $I \subset E$ be a two-sided ideal. Let $J \subset R$ be the ideal generated by the non-constant coefficients of the polynomials in the set
\[
\{D_{R[t]}(1-xt) \ | \ x \in I\} \subset R[t].
\]
Let $E' := E/(I,J)$ and $R' := R/J$, and let $f:R \to R'$ and $g:E \to E'$ be the quotient maps. 
\begin{enumerate}
\item There exists a unique Cayley-Hamilton pseudorepresentation $D' : E' \ra R'$ such that the pair $(f,g)$ gives a morphism $(E,D) \to (E',D')$ of Cayley-Hamilton algebras.
\item For any Cayley-Hamilton representation $\rho: (E,D) \ra (E_A, D_A : E_A \ra A)$, the map $E \to E_A$ sends $I$ to $0$ if and only if $\rho$ factors through $(E,D) \ra (E',D')$. 
\end{enumerate}
\end{lem}
\begin{proof}
(1)\ The uniqueness follows from (2). To show existence, we start with the pseudorepresentation $\tilde{D}:=(f \circ D): E\otimes_R R' \to R'$. To construct $D'$, we use Chenevier's notion of \emph{kernel of a pseudorepresentation} \cite[\S 1.17]{chen2014}. It is defined by the following universal property: there is an $R'$-algebra quotient $h:E \otimes_R R' \to (E \otimes_R R')/\ker(\tilde{D})$ and a pseudorepresentation $\tilde{D}': (E \otimes_R R')/\ker(\tilde{D}) \to R'$ such that $\tilde{D}=\tilde{D}' \circ h$, and $\ker(\tilde{D})$ is the maximal ideal with this property.

There is an equality
\begin{equation}
\ker(\tilde{D}) = \{x \in E\otimes_R R' \ | \ \tilde{D}_{R'[t]}(1-xyt)=1 \ \forall y \in E \otimes_R R'\}.
\end{equation}
This is proven in \cite[Lem.\ 2.1.2]{CS2016}\footnote{This lemma was removed from the final version \cite{CS2019}, having become unnecessary.}. Their proof uses Amitsur's formula \cite{amitsur1980} for pseudorepresentations \cite[(1.5)]{chen2014}, together with the fact that, for any $x,y \in E$, $D_{R[t]}(1-xyt)=D_{R[t]}(1-yxt)$, which can be proven in the same way as for usual determinants in linear algebra.

By the definition of $J$, we have $\tilde{D}_{R'[t]}(1-(x\otimes 1)t)=1$ for all $x \in I$, and hence $\tilde{D}_{R'[t]}(1-(x \otimes 1)yt)=1$ for all $x \in I$ and $y \in E\otimes_R R'$, since $I$ is a two-sided ideal. This implies that $I\otimes_R R' \subset \ker(\tilde{D})$. By the property of $\ker(\tilde{D})$, we have $D':E' \to R'$ such that $f \circ D = D' \circ g$. 

(2)\ If $\rho$ factors through $(E,D) \ra (E',D')$, then the map $E \to E_A$ factors through $E \to E' \to E_A$, so $E \to E_A$ sends $I$ to $0$. Conversely, let $\rho:(E,D) \ra (E_A,D_A)$ be given by $f_A:R \to A$ and $g_A:E \to E_A$, and assume that $g_A(I)=0$. To show that $\rho$ factors through $(E,D) \ra (E',D')$, it suffices to show that $f_A(J)=0$. For any $x \in E$, the naturality of $D$ implies that the image of $D_{R[t]}(1-xt)$ in $A[t]$ under $f_A$ is given by $D_{A[t]}(1-xt) = (f_A \circ D)_{A[t]}(1-xt)$. Hence it is enough to show that, for $x \in I$, we have $(f_A \circ D)_{A[t]}(1-xt)=1$. However, since $\rho$ is a morphism of Cayley-Hamilton algebras, we have
\[
(f_A \circ D)_{A[t]}(1-xt) = (D_A \circ g_A)_{A[t]}(1-xt) = (D_A)_{A[t]}(g_A(1-xt)) = (D_A)_{A[t]}(1)=1,
\]
where we use the fact that $g_A(x)=0$. 
\end{proof}

\begin{eg}
\label{eg:CH-quot-with-scalar-ideal}
Let $(E,D)$ be a Cayley-Hamilton $R$-algebra and let $I=JE$ with $J \subset R$ being an ideal. Then the Cayley-Hamilton quotient $(E',D')$ of $(E,D)$ by $I$ has $E'=E/JE$ and scalar ring $R'=R/J$. In other words, $(E', D') = (E/JE, D \otimes_R R/J : E/JE \ra R/J)$. 
\end{eg}
\begin{defn}
\label{defn:ERP}
With the notation of the lemma, we call $(E',D')$ the \emph{Cayley-Hamilton quotient of $(E,D)$ by $I$}. 
\end{defn}

\begin{defn}
\label{defn:E^C_Db}
Let $K^\cC \subset E_\Db$ be the kernel of the algebra homomorphism \eqref{eq:KP}. Let $(E^\cC_\Db, D_{E^\cC_\Db})$ denote the Cayley-Hamilton quotient of $(E_\Db, D^u_{E_\Db})$ by $K^\cC$, and let $R_\Db^\cC$ denote the scalar ring of $E^\cC_\Db$. 
\end{defn}

\subsection{Extending condition $\cC$ to pseudorepresentations and Cayley-Hamilton representations}
\label{subsec:CH-P}

We extend $\cC$ to $A[G]$-modules that may not have finite cardinality in the following way. 

\begin{defn}
\label{defn:proP}
Let $\cC$ be a stable condition on objects of $\Mt$.  Let $(A,\m_A)$ be a complete commutative Noetherian local $\bZ_p$-algebra. For an $A[G]$-module $M$ that is finitely generated as an $A$-module, we say that $M$ \emph{satisfies condition $\cC$} when $M/\m_A^i M$ satisfies $\cC$ for all $i \geq 1$. 
\end{defn}

Note that, for $(A,\m_A)$ and $M$ as in the definition, the canonical map $M \to \varprojlim_i M/\m_A^i M$ is an isomorphism. Using this, one can check that this extension of $\cC$ is stable in the same sense as Definition \ref{defn:stable}. We will use this extension of condition $\cC$ without further comment. 

Now we give the definition of condition $\cC$ for Cayley-Hamilton representations. 
\begin{defn}
\label{defn:PforCH}
	Let $(A,\m_A)$ be a complete commutative Noetherian local $\bZ_p$-algebra and let $(E, \rho,D)$ be a Cayley-Hamilton representation of $G$ with scalar ring $A$ and residual pseudorepresentation $\Db$. We say that $(E, \rho, D)$ \emph{satisfies condition $\cC$} if $E$ satisfies condition $\cC$ as an $A[G]$-module. (Note that $E$ is finitely generated as an $A$-module by Proposition \ref{prop:C-H alg is module-finite}.)
	
	We let $\CH_{G,\Db}^\cC$ denote the full subcategory of $\CH_{G,\Db}$ whose objects satisfy condition $\cC$.
\end{defn}

We can be confident that this notion behaves well by finding a universal object. 
\begin{thm}
\label{thm:univ_CH_P}
The Cayley-Hamilton representation $(E^\cC_\Db, D_{E^\cC_\Db} : E^\cC_\Db \ra R^\cC_\Db)$ of Definition \ref{defn:E^C_Db} is the universal object in $\CH_{G,\Db}^\cC$.
\end{thm}

\begin{proof}
By the definition of $(E^\cC_\Db, D_{E^\cC_\Db} : E^\cC_\Db \ra R^\cC_\Db)$, we see that the map $E_\Db \to E_\Db^\cC$ sends $K^\cC$ to $0$. By Lemma \ref{lem:EP}, we see that $(E^\cC_\Db, D_{E^\cC_\Db} : E^\cC_\Db \ra R^\cC_\Db)$ satisfies condition $\cC$. We now show that it has the universal property.

Let $A$ be a complete commutative Noetherian local $\bZ_p$-algebra, and let $(E, \rho, D)$ be a Cayley-Hamilton representation with scalar ring $A$ and residual pseudorepresentation $\Db$. We have to show that $(E, \rho, D)$ satisfies condition $\cC$ if and only if the map of Cayley-Hamilton algebras $(E_\Db, D^u_{E_\Db}) \ra (E,  D)$ induced by the universal property of $(E_\Db, \rho^u, D^u_{E_\Db})$ factors through $(E^\cC_\Db, D_{E^\cC_\Db})$. 

	For any $i \geq 1$, $E_\Db \ra E/\m_A^i E$ factors through $E_\Db/\mDb^i E_\Db$. (Recall that a local homomorphism of scalar rings $R_\Db \ra A$ is implicit in $(E_\Db, D^u_{E_\Db}) \ra (E,  D)$.) 
		
	By Lemma \ref{lem:EP}, $(E, \rho, D)$ satisfies $\cC$ if and only if $E_\Db/\mDb^{i} E_\Db \ra E/\m_A^i E$ factors through $E^\cC_\Db(i)$ for every $i \geq 1$. Equivalently, $(E, \rho, D)$ satisfies $\cC$ if and only if $E_\Db \to E$ maps $K^\cC$ to $0$. By Lemma \ref{lem:max_CH_quot},  $K^\cC$ maps to $0$ in $E_\Db \ra E$ if and only if $(E_\Db, D^u_{E_\Db}) \ra (E, D)$ factors through $(E^\cC_\Db, D_{E^\cC_\Db})$. 
\end{proof}

Following the pattern of \cite[Defn.\ 5.9.1]{WWE1}, we define condition $\cC$ on pseudorepresentations. 

\begin{defn}
\label{defn:PsR-P}
	Let $A$ be an complete commutative Noetherian local $\bZ_p$-algebra. Let $D : G \ra A$ be a pseudorepresentation with residual pseudorepresentation $\Db$. Then $D$ \emph{satisfies condition $\cC$} provided that there exists a Cayley-Hamilton representation $(E, \rho, D')$ satisfying condition $\cC$ such that $D = \psi(\rho) := D' \circ \rho$. 
\end{defn}

We define the $\cC$-pseudodeformation functor $\PsDef^\cC_\Db: \hat\cC_{W(\F)} \to \mathrm{Sets}$ by 
\[
\PsDef^\cC_\Db(A)= \{\text{pseudodeformations } D : G \ra A \text{ of }\Db \text{ satisfying } \cC\}. 
\]

\begin{thm}
\label{thm:RDP}
The functor $\PsDef^\cC_\Db$ is represented by $R^\cC_\Db$. 
\end{thm}

\begin{proof}
Let $A \in \hat\cC_{W(\F)}$, and let $D \in \PsDef_\Db(A)$. By Theorem \ref{thm:RDb}, there is unique $R_\Db \ra A$ such that $D \cong D^u_\Db \otimes_{R_\Db} A$. We have to show that $D \in \PsDef^\cC_\Db(A)$ if and only if $R_\Db \ra A$ factors through $R_\Db \onto R_\Db^\cC$.

If $R_\Db \ra A$ factors through $R_\Db \rsurj R^\cC_\Db$, then the Cayley-Hamilton representation 
\[
(E^\cC_\Db \otimes_{R^\cC_\Db} A, \rho^\cC \otimes_{R_\Db} A : G \ra (E^\cC_\Db \otimes_{R^\cC_\Db} A)^\times, D_{E^\cC_\Db} \otimes_{R^\cC_\Db} A)
\]
induces $D$ via $D = (R^\cC_\Db \ra A) \circ D_{E^\cC_\Db}$ and satisfies condition $\cC$ by Theorem \ref{thm:univ_CH_P}. Consequently $D$ satisfies $\cC$.  

Now assume $D$ satisfies condition $\cC$, i.e.\ there exists a Cayley-Hamilton representation $(E, \rho, D')$ satisfying $\cC$ such that $D = D' \circ \rho$. By Theorem \ref{thm:univ_CH_P}, there exists a morphism of Cayley-Hamilton algebras $(E^\cC_\Db, D_{E^\cC_\Db}) \ra (E, D)$ inducing $\rho$. In particular, the implicit morphism of scalar rings $R^\cC_\Db \ra A$ factors $R_\Db \ra A$.
\end{proof}

\subsection{Modules with Cayley-Hamilton structure}
\label{subsec:ModCH-P}
We introduce the notion of Cayley-Hamilton $G$-module. 
\begin{defn} 
Let $A \in \hat \cC_{W(\bF)}$. A \emph{Cayley-Hamilton $G$-module} over $A$ is the data of a Cayley-Hamilton representation $(E,\rho,D)$ of $G$ with scalar ring $A$, and an $E$-module $N$. We consider $N$ as a $A[G]$-module via the map $\rho: A[G] \to E$. We often refer to a to a Cayley-Hamilton $G$-module simply by the letter $N$, and call $(E,D)$ the \emph{Cayley-Hamilton algebra of $N$}. We say $N$ is \emph{faithful} if it is faithful as $E$-module. 
\end{defn}

\begin{eg}
If $N$ is an $A[G]$-module and there is a Cayley-Hamilton pseudorepresentation $D: \End_A(N) \to A$, then the canonical action of $\End_A(N)$ on $N$ makes $N$ a faithful Cayley-Hamilton $G$-module with Cayley-Hamilton algebra $(\End_A(N),\rho,D)$, where $\rho:G \to \End_A(N)$ is the action map.
\end{eg}

\begin{eg}
\label{eg:Tate module}
As a special case of the previous example, suppose $N=N_1 \oplus N_2$ as $A$-modules such that $\End_A(N_i)=A$ (i.e.~the only endomorphisms are scalars). Note that $N_1$ and $N_2$ need not be free $A$-modules (for example, they could be dualizing $A$-modules with $A$ non-Gorenstein; or, if $A=\Z_p\lb T\rb$, they could be the maximal ideal).
Then $\End_A(N)$ has the structure of $A$-GMA $(B,C,m)$ as in Example \ref{eg:GMA1,1}, where $B=\Hom_A(N_1,N_2)$ and $C=\Hom_A(N_2,N_1)$, and $m(f,g)=f\circ g$ (see \S\ref{sec:GMAs}, below, for a discussion of GMAs). In particular, there is a natural Cayley-Hamilton pseudorepresentation $D: \End_A(N) \to A$.

This example appears frequently in the study of pseudorepresentations associated to ordinary modular forms, where one takes $A$ to be a Hecke algebra, and $N$ to be the $p$-adic Tate module of a modular Jacobian, in which case $N$ is a direct sum of a free $A$-module of rank $1$ and a dualizing $A$-module (c.f.~\cite{mazur1978,WWE1,WWE2,WWE3}).
\end{eg}

\begin{thm}
\label{thm:ModToCH}
Let $N$ be a faithful Cayley-Hamilton $G$-module with Cayley-Hamilton $A$-algebra $(E,D)$. Then $N$ satisfies condition $\cC$ as an $A[G]$-module if and only if $(E,\rho,D)$ satisfies condition $\cC$ as a Cayley-Hamilton representation.
\end{thm}
\begin{proof}
By Definition \ref{defn:proP}, it suffices to prove the theorem in the case that $A$ is Artinian and local. Choose $i \geq 1$ such that $\mDb^i \cdot A = 0$. Let $\Db:G \to \bF$ be the residual pseudorepresentation of $D \circ \rho: G \ra A$.

By Theorem \ref{thm:univ_CH}, there is a distinguished morphism of Cayley-Hamilton algebras $(E_\Db, D^u_{E_\Db}) \ra (E,D)$. Then the action map $\Z_p[G] \to \End_A(N)$ factors as 
\begin{equation}
\label{eq:factoring-G-action}
\Z_p[G] \to E_\Db \to E_\Db/\m_\Db^iE_\Db \to E \to \End_A(N).
\end{equation}
We claim that the following are equivalent: 
\begin{enumerate}
\item The $\Z_p[G]$-module $N$  satisfies $\cC$
\item The map $E_\Db/\m_\Db^iE_\Db \to \End_A(N)$ factors through $E_\Db^\cC(i)$
\item The map $E_\Db  \to \End_A(N)$ sends $K^\cC$ to $0$
\item The map $E_\Db  \to E$ sends $K^\cC$ to $0$
\item The morphism of Cayley-Hamilton algebras $(E_\Db, D^u_{E_\Db}) \ra (E,D)$ factors through $(E_\Db^\cC,D_{E^\cC_\Db})$
\item The Cayley-Hamilton representation $(E,D)$ satisfies condition $\cC$.
\end{enumerate}
(Recall that $K^\cC$ was defined in Definition \ref{defn:E^C_Db}.) The equivalences are proven as follows:
\begin{itemize}[leftmargin=1.75em]
\item[]$(1) \Longleftrightarrow (2)$: Lemma \ref{lem:EP}(3).
\item[]$(2) \Longleftrightarrow (3)$: From the definition of $K^\cC$ and \eqref{eq:factoring-G-action}.
\item[]$(3) \Longleftrightarrow (4)$: Since $N$ is faithful, the map $E \to \End_A(N)$ is injective.
\item[]$(4) \Longleftrightarrow (5)$: Lemma \ref{lem:max_CH_quot}.
\item[]$(5) \Longleftrightarrow (6)$: Theorem \ref{thm:univ_CH_P}. \qedhere
\end{itemize}
\end{proof}

\subsection{Formal moduli of representations with $\cC$}
\label{subsec:formal moduli with C}
We generalize Ramakrishna's result (Theorem \ref{thm:Ram}) to any family of integral $p$-adic representations of $G$. 

\begin{thm}
\label{cor:P-closed}
	There exists a unique closed formal substack (resp.\ closed formal subscheme)
	\[
 \Rep^{\cC,d}_G \subset \Rep_G^d	, \qquad (\text{resp.\ }\Rep^{\square,\cC,d}_G \subset \Rep^{\square,d}_G)
	\]
	characterized by the following property. For any commutative local $\bZ_p$-algebra $B$ of finite cardinality and free rank $d$ $B$-module $V_B$ with a $B$-linear $G$-action (resp.\ and fixed basis), the corresponding $B$-point of $\Rep_G^d$ lies in $\Rep^{\cC,d}_G$ (resp.\ of $\Rep_G^{\square,d}$ lies in $\Rep^{\square,\cC,d}_G$) if and only if $V_B$ has $\cC$. 
\end{thm}
\begin{proof}
It suffices to consider the case of $\Rep^{\square,d}_G$. Indeed, since condition $\cC$ does not depend on the choice of basis, the closed subscheme $\Rep^{\square,\cC,d}_G \subset \Rep^{\square,d}_G$ descends to a closed locus in $\Rep_G^d$. 

By Theorem \ref{thm:RDb}, we may consider a fixed residual pseudorepresentation $\Db$ and produce $\Rep^{\square,\cC}_\Db \subset \Rep^\square_\Db$. We define $\Rep^{\square,\cC}_\Db$ via the pullback diagram
\[\xymatrix{
\Rep^{\square,\cC}_\Db \ar[r] \ar[d] & \Rep^\square_\Db \ar[d]^-\wr \\
\Rep^\square_{E_\Db^\cC,D_{E^\cC_\Db}} \ar[r] & \Rep^\square_{E_\Db,D^u_{E_\Db}}
}\]
where the right vertical arrow is the isomorphism in Theorem \ref{thm:ED-compat}. 

Let $B$ be an local $\bZ_p$-algebra of finite cardinality. By definition, a point $V_B \in \Rep^\square_\Db(B)$ lies in $\Rep^{\square,\cC}_\Db(B)$ if and only if the map $(E_\Db,D^u_{E_\Db}) \to (\End_B(V_B),\det)$ factors through $(E_\Db^\cC,D_{E^\cC_\Db})$. By Theorem \ref{thm:univ_CH_P}, this occurs if and only if the Cayley-Hamilton representation $(\End_B(V_B),\det)$ has $\cC$, which, by Theorem \ref{thm:ModToCH}, is if and only if $V_B$ has $\cC$. By Lemma \ref{lem:tfg algebra}, this characterizes $\Rep^{\square,\cC}_\Db$.
\end{proof}

\begin{cor}
\label{thm:univ_EP}
Let $A\in \tfg_{\Z_p}$ and $(V_A,\rho_A) \in \Rep^d_G(A)$. There exists a unique quotient morphism $A \rsurj A^\cC$ in $\tfg_{\Z_p}$ such that, for any local $A$-algebra $B$ of finite cardinality, the object $V_A \otimes_A B$ of $\Mt$ satisfies $\cC$ if and only if the homomorphism $A \ra B$ factors through $A^\cC$.
\end{cor}

\begin{proof}
The quotient $A \to A^\cC$ is defined by the pullback square
\[
\xymatrix{
\Spec(A^\cC) \ar[r] \ar[d] & \Rep^{\cC,d}_G \ar[d] \\
\Spec(A) \ar[r] & \Rep^{d}_G.
}
\]
The characterizing property follows from the characterizing property of $\Rep^{\cC,d}_G$.
\end{proof}

\section{Deformation conditions for generalized matrix algebras} 
\label{sec:GMAs}

For the remainder of the paper, we assume that the pseudorepresentation $\Db$ is multiplicity-free (see Definition \ref{defn:RMF}). Under this assumption, as was first noticed by Chenevier \cite{chen2014}, the universal Cayley-Hamilton algebra with residual pseudorepresentation $\Db$ admits the additional structure of a generalized matrix algebra (GMA). Generalized matrix algebras, which were first introduced by Bella\"iche and Chenevier \cite{BC2009}, are a particularly concrete and explicit type of Cayley-Hamilton algebra. 

In this section, we recall in \S \ref{subsec:GMAs} the theory of generalized matrix algebras and in \S \ref{subsec:RMF} why the universal Cayley-Hamilton algebra admits this structure. In \S\ref{subsec:RMF-P}, we exploit this extra structure to prove that the universal Cayley-Hamilton algebra with property $\cC$, constructed in \S \ref{sec:prop}, is characterized by representations. This is the main new result of the section. Explicitly, in the multiplicity-free case, we show that, if every representation with $\cC$ factors through a quotient of universal Cayley-Hamilton algebra, then that quotient must be universal Cayley-Hamilton algebra with property $\cC$. 

The results of \S \ref{subsec:GMAs} and \S \ref{subsec:RMF} will also be used in \S\ref{sec:gma_ext}, where we will exploit the GMA structure to compute the Cayley-Hamilton algebra with property $\cC$ in terms of group cohomology.

\subsection{Generalized matrix algebras and their adapted representations}
\label{subsec:GMAs}

A generalized matrix algebra is a particular kind of Cayley-Hamilton algebra with extra data. We learned this notion from Bella\"iche--Chenevier \cite{BC2009}. 
\begin{defn}[{\cite[\S1.3]{BC2009}}]
\label{defn:GMA}
	Let $R$ be a ring. A \emph{generalized matrix algebra over $R$} (or an \emph{$R$-GMA}) is the data of:
	\begin{enumerate}
		\item An $R$-algebra $E$ that is finitely generated as an $R$-module,
		\item A set of orthogonal idempotents $e_1, \dotsc, e_r \in E$ such that $\sum_i e_i =1$, and
		\item A set of isomorphisms of $R$-algebras $\phi_i : e_i E e_i \risom M_{d_i}(R)$ for $i=1,\dots,r$.
	\end{enumerate}
	We call $\cE = (\{e_i\}, \{\phi_i\})$ the \emph{data of idempotents} or \emph{GMA structure} of $E$, and write the $R$-GMA as $(E,\cE)$.  We call the list of integers $(d_1, \dots, d_r)$ the \emph{type} of $(E,\cE)$. These data are required to satisfy the condition that the function $\Tr_\cE:E \ra A$ defined by 
	\[
	\Tr_\cE(x) := \sum_{i=1}^r \tr (\phi_i(e_ixe_i))
	\]
	is a central function, i.e.\ $\Tr_\cE(xy) = \Tr_\cE(yx)$ for all $x,y \in E$. 
	
	Given an $R$-GMA $(E,\cE)$ and an $R'$-GMA $(E', \cE')$, a \emph{morphism of GMAs} $\rho:(E,\cE) \ra (E', \cE')$ is the data of a ring homomorphism $f:R \to R'$ and morphism of $R$-algebras $\rho: E \ra E'$ such that $\cE$ and $\cE'$ are of the same type $(d_1,\dots,d_r)$, we have $g(e_i)=e_i'$ and $f \circ \phi_i = \phi_i'\circ g$ for $i=1,\dots,r$. 
\end{defn}

\begin{eg}
\label{eg:matrix}
The matrix algebra $M_d(R)$ comes with a natural $R$-GMA structure $\cE = (1, \mathrm{id}: M_d(R) \risom M_d(R))$ of type $(d)$. More generally, given any ordered partition of $d = d_1 + \dotsc + d_r$ of $d$, there is a natural $R$-GMA structure $\cE_{\mathrm{block}}$ on $M_d(R)$ of type $(d_1, \dotsc, d_r)$. Namely, the natural $R$-algebra map with block-diagonal image 
\[
\nu_1 \times \dotsm \times \nu_r : M_{d_1}(R) \times \dotsm \times M_{d_r}(R) \rinj M_d(R)
\]
induces $\cE_{\mathrm{block}} = (e_i = \nu_i(1_i), \phi_i)$, where $1_i \in M_{d_i}(R)$ is the identity matrix, and $\phi_i$ is the inverse to $\nu_i$ on its image $e_iM_d(R)e_i$. 
\end{eg}

\begin{lem}
\label{lem:DE}
Given an $R$-GMA $(E,\cE)$, there is a canonical Cayley-Hamilton pseudorepresentation $D_\cE: E \ra R$, such that $\Tr_{D_\cE}=\Tr_\cE$. A morphism of $R$-GMAs $(E,\cE) \to (E',\cE')$ induces a morphism of Cayley-Hamilton algebras $(E,D_\cE) \to (E',D_{\cE'})$.
\end{lem} 
\begin{proof}
See \cite[Prop.\ 2.23]{WE2018}; its statement includes the first claim. The second claim follows from examining the formula for $D_\cE$ given in \emph{loc.\ cit.}, noting that a morphism of GMAs preserves the idempotents that are used to specify $D_\cE$.
\end{proof}

This lemma gives a faithful embedding of the category of $R$-GMAs into the category of $R$-Cayley-Hamilton algebras. We will consider this embedding as an inclusion. We extend the definition of Cayley-Hamilton representation to GMAs as follows.
\begin{defn}
If $(E',\cE')$ is a GMA, we refer to a morphism of pseudorepresentations $(E,D) \to (E',D_{\cE'})$ as a \emph{GMA representation} of $(E,D)$. If $(E,\cE)$ is another GMA, we call a GMA representation $(E,D_\cE) \to (E',D_{\cE'})$ \emph{adapted} if the same data give a morphism of GMAs $(E, \cE) \to (E',\cE')$.
\end{defn}

\begin{lem}
\label{lem:gma data}
Given an $R$-GMA $(E,\cE)$ of type $(d_1,\dots,d_r)$, we can associate to it the data of
\begin{enumerate}
\item $R$-modules $\cA_{i,j}$ for $1 \le i,j \le r$,
\item canonical isomorphisms $\cA_{i,i} \isoto R$ for $1 \le i \le r$, and
\item $R$-module homomorphisms $\varphi_{i,j,k}:\cA_{i,j} \otimes_R \cA_{j,k} \ra \cA_{i,k}$ for $1 \leq i,j,k \leq r$,
\end{enumerate}
such that $(\cA_{i,j},\varphi_{i,j,k})$ completely determine $(E,\cE)$ and there is an isomorphism of $R$-modules
\begin{equation}
	\label{eq:GMA_form}
	E \lrisom \begin{pmatrix}
	M_{d_1}(\cA_{1,1}) & M_{d_1 \times d_2}(\cA_{1,2}) & \dotsm & M_{d_1 \times d_r}(\cA_{1,r}) \\
	M_{d_2 \times d_1}(\cA_{2,1}) & M_{d_2}(\cA_{2,2}) & \dotsm & M_{d_2 \times d_r}(\cA_{2,r}) \\
	\vdots & \vdots & \vdots & \vdots \\
	M_{d_r \times d_1}(\cA_{r,1}) & M_{d_2}(\cA_{r,2}) & \dotsm & M_{d_r}(\cA_{r,r}) \\
	\end{pmatrix},
	\end{equation}

Moreover, the collection of maps $\varphi_{i,j,k}$ satisfies properties (UNIT), (COM), and (ASSO) given in \cite[\S1.3.2]{BC2009}, and there is a bijection between $R$-GMAs of type $(d_1,\dots,d_r)$
and data $(\cA_{i,j},\varphi_{i,j,k})$ satisfying (UNIT), (COM), and (ASSO).
\end{lem}
\begin{proof}
This is explained in \cite[\S1.3.1-\S1.3.6]{BC2009}. The association is given as follows. Let $E_i:= \phi_i^{-1}(\delta^{1,1})$, where $\delta^{1,1}$ denotes the elementary matrix with 1 as $(1,1)$th entry, and 0 otherwise. Define $\cA_{i,j} := E_j E E_i$. The maps $\varphi_{i,j,k}$ are induced by the multiplication in $E$. In particular, note that $\phi_i$ induces a canonical isomorphism $\cA_{i,i} \isoto R$.
\end{proof}

We will not spell out the bijection or the properties (UNIT), (COM), and (ASSO) in general. Instead, we explain the content of the lemma in the case of type $(1,1)$.

\begin{eg}
\label{eg:GMA1,1}
There is a bijection between $R$-GMAs $(E,\cE)$ of type $(1,1)$ and triples $(B,C, m)$ where $B,C$ are finitely generated $R$-modules and $m : B \otimes_R C \ra R$ is an $R$-module homomorphism, such that the squares
\[
\xymatrix{
B \otimes_R C \otimes_R B \ar[r]^-{1 \otimes (m \circ \iota)} \ar[d]^-{m \otimes 1} & B \otimes_R R \ar[d] & & 
C \otimes_R B \otimes_R C \ar[r]^-{1 \otimes m} \ar[d]^-{(m \circ \iota) \otimes 1} & C \otimes_R R \ar[d]  \\ 
R \otimes_R B \ar[r] & B & & R \otimes_R C \ar[r] & C
}
\]
commute. Here $\iota: C \otimes_R B \isoto B \otimes_R C$ is the isomorphism given by $b \otimes c \mapsto c \otimes b$, and the unlabeled maps are the $R$-action maps.

The $R$-GMA associated to a triple $(B,C,m)$ is 
\begin{equation}
\label{eq:E-form}
E =\begin{pmatrix} 
R & B \\ C & R
\end{pmatrix}.
\end{equation}
This means that $E=R \oplus B \oplus C \oplus R$ as an $R$-module, and the multiplication on $E$ is given by $2\times2$-matrix multiplication, but where the action maps $R \otimes_R B \to B$, $R \otimes_R C \to C$ and the map $m$ are used in place of the scalar multiplication. The idempotents $e_1, e_2$ are given in this notation by $e_1=\sm{1}{0}{0}{0}$ and $e_2=\sm{0}{0}{0}{1}$, and the isomorphisms $\phi_i$ are given by the identifications $\sm{R}{0}{0}{0} \isoto R$ and $\sm{0}{0}{0}{R} \isoto R$.

In the notation of the lemma, $\cA_{1,2}=B$ and $\cA_{2,1}=C$, and the maps $\varphi_{i,j,k}$ are the action maps, except $\varphi_{1,2,1}=m$ and $\varphi_{2,1,2}= m \circ \iota$. Indeed, the (UNIT) property requires that the maps $\varphi_{i,j,k}$ equal the action maps, unless $(i,j,k) \in \{(1,2,1),(2,1,2)\}$. The (ASSO) property is expressed by the commutative squares above, and the (COM) property is that 
$\varphi_{2,1,2}= \varphi_{1,2,1} \circ \iota$.
\end{eg}

\begin{eg}
\label{eg:m0}
In the foregoing example, setting $m$ to be the zero map is always a valid choice. 
\end{eg}

\begin{defn}[{\cite[Defn.\ 1.3.6]{BC2009}}]
	\label{defn:adapted}
	Let $A$ be a commutative $R$-algebra and let $(E,\cE)$ be an $R$-GMA of type $(d_1, \dotsc, d_r)$. Let $d = \sum_{i=1}^r d_i$ and let $(M_d(A),\cE_\mathrm{block})$ be the $A$-GMA of type $(d_1, \dotsc, d_r)$ constructed in Example \ref{eg:matrix}. 
	
	An \emph{adapted representation of $(E,\cE)$}, denoted $\rho:(E,\cE) \to M_d(A)$ is a morphism of GMAs $\rho:(E,\cE) \to (M_d(A),\cE_\mathrm{block})$. A pseudorepresentation $D :E \otimes_R A \ra A$ is \emph{adapted to $\cE$} if $D=D_\cE \otimes_R A$. When the data of idempotents $\cE$ is understood, we say $D$ is \emph{adapted} instead of adapted to $\cE$. 
\end{defn}

Now fix $(E,\cE)=(E,(\{e_i\},\{\phi_i\}))$, an $R$-GMA of type $(d_1, \dotsc, d_r)$, and let $(\cA_{i,j},\varphi_{i,j,k})$ be the data associated to it by Lemma \ref{lem:gma data}. We define a the set-valued functor on commutative $R$-algebras by
	\[
	\Rep_{(E,\cE),\Ad}^\square:A \mapsto \{\text{Adapted representations } \rho:(E,\cE) \to M_d(A)\}.
	\]
There is an explicit presentation for $\Rep^\square_{(E,\cE),\Ad}$ as an affine $\Spec R$-scheme as follows. 

\begin{thm}
\label{thm:Ad_main}
The functor $\Rep_{(E,\cE),\Ad}^\square$ is represented by $\Spec S$, where $S$ is an $R$-algebra quotient 
\[
\Sym^*_R \left(\bigoplus_{1 \leq i \neq j \leq r} \cA_{i,j} \right) \rsurj S
\]
by the ideal generated by the set of all $\varphi(b \otimes c)-b \otimes c$, with $b \in \cA_{i,j}, c \in \cA_{j,k}$ and $\varphi=\varphi_{i,j,k}$ for all $1 \le i, j, k \leq r$. In particular, $S$ is a finitely generated $R$-algebra.

Moreover, for $1 \leq i, j \leq r$, the natural $R$-module maps $\cA_{i,j} \ra S$ are split injections of $R$-modules, inclusive of the case $R \risom \cA_{i,i} \rinj S$, which is the $R$-algebra structure map of $S$. The universal adapted representation $\rho_\Ad:(E,\cE) \ra M_d(S)$ is given by the isomorphism of \eqref{eq:GMA_form} along with these $R$-module injections. In particular, the $R$-algebra homomorphism $E \to M_d(S)$ is injective. 
\end{thm}

\begin{proof}
See \cite[Prop.\ 1.3.9]{BC2009} and its proof, as well as  \cite[Prop.\ 1.3.13]{BC2009} for the split injectivity. 
\end{proof}

By \cite[Prop.\ 2.23]{WE2018}, any adapted representation of $(E,\cE)$ is compatible with $D_\cE$. This determines a monomorphism $\Rep_{(E,\cE),\Ad}^\square\rinj \Rep^\square_{E,D_\cE}$, which can easily be observed to be a closed immersion. 

Let $\rho^\square$ denote the universal object of $\Rep^\square_{E,D_\cE}$. Let $Z(\cE) \subset \GL_d$ denote the subgroup that stabilizes $\rho^\square(e_iEe_i)$, for all $i=1,\dots,r$, under the adjoint action of $\GL_d$ on $\Rep^\square_{E,D_\cE}$. This $Z(\cE)$ is a torus of rank $r$. 

\begin{prop}[{\cite[Thm.\ 2.27]{WE2018}}]
	\label{prop:adapted_GMA}
For any $R$-GMA $(E,\cE)$, the map $\Rep_{(E,\cE),\Ad}^\square\rinj \Rep^\square_{E,D_\cE}$ induces an isomorphism	
\begin{equation}
\label{eq:adapted_to_reg}
	[\Rep_{(E,\cE),\Ad}^\square / Z(\cE)] \lrisom \Rep_{E,{D_\cE}} 
\end{equation}
of $\Spec R$-algebraic stacks. 
\end{prop}

\subsection{Residually multiplicity-free representations of profinite groups}
\label{subsec:RMF}

Let $G$ be a profinite group satisfying the $\Phi_p$ finiteness condition, and let $\bF$ be a finite field of characteristic $p$ with algebraic closure $\overline{\bF}$. 

By \cite[Cor.\ 2.9(2)]{WE2018}, which is a mild strengthening of \cite[Thm.~2.16]{chen2014}, for any pseudorepresentation $\Db: G \ra \bF$, there is a unique semi-simple representation $\rho^{ss}_\Db :G \to \GL_d(\bF)$ such that $\Db =\det \circ \rho^{ss}_\Db$. 

\begin{defn}
\label{defn:RMF}
	A residual pseudorepresentation $\Db: G \ra \bF$ is \emph{multiplicity-free} if $\rho^{ss}_\Db \otimes_\F \overline{\F}$ has pairwise non-isomorphic simple factors, and each of these factors is defined over $\F$. Equivalently, $\rho^{ss}_\Db$ is a direct sum of pairwise non-isomorphic absolutely irreducible representations. A representation $(V_A,A) \in \Rep_G^d(A)$ is \emph{residually multiplicity-free} if $(V_A,A) \in \Rep_\Db(A)$ with $\Db$ multiplicity-free.
\end{defn}

Note that when $\rho^{ss}_\Db \otimes_\F \overline{\F}$ has distinct simple factors, then $\rho^{ss}_\Db$ is multiplicity-free after replacing $\F$ with a finite extension.

\begin{thm}
\label{thm:CH-GMA}
	Let $\Db: G \to \bF$ be multiplicity-free, and let $(d_1,\dots,d_r)$ be the dimensions of the simple factors of $\rho^{ss}_\Db$. Let $A$ be a Noetherian Henselian local ring with residue field $\bF$, and let $(E,\rho,D)$ Cayley-Hamilton representation over $A$ with residual pseudorepresentation $\Db$. Then there is an $A$-GMA structure $\cE$ of type $(d_1,\dots,d_r)$ on $E$ such that $D=D_{\cE}$. 
	
Moreover, given a morphism $(E,\rho,D) \to (E',\rho',D')$ of such objects, the structures $\cE$ and $\cE'$ may be chosen so that the map $(E,\cE) \to (E',\cE')$ is a morphism of GMAs.
\end{thm}
\begin{proof}
The structure $\cE$ is constructed in \cite[Thm.\ 2.22(ii)]{chen2014}, and it follows from this construction that $D= D_{\cE}$ (see \cite[Thm.\ 2.27]{WE2018}). Moreover, the construction only depends on the choice certain lifts of idempotents. If we first choose the structure $\cE=(\{e_i\},\{\phi_i\})$ on $E$, then the images of $e_i$ in $E'$ will give a choice of lifts of idempotents in $E'$, and for the resulting GMA structure $\cE'$, the map $(E,\cE) \to (E',\cE')$ is a morphism of GMAs.
\end{proof}

For the rest of this section, we fix a pseudorepresentation $\Db: G \to \bF$ that is multiplicity-free. By Theorem \ref{thm:RDb}, the ring $R_\Db$ is Noetherian and complete (and hence Henselian). By Theorem \ref{thm:CH-GMA}, we can and do fix a choice of $A$-GMA structure $\cE_\Db$ of type $(d_1,\dots,d_r)$ on $E_\Db$ such that $D^u_{E_\Db}=D_{\cE_\Db}$.

\begin{cor}
\label{cor:compat=adapt}
Assume that $\Db$ is multiplicity-free, and let $\cE_\Db$ be choice of $R_\Db$-GMA structure on $E_\Db$ as in Theorem \ref{thm:CH-GMA}.
	
	\begin{enumerate}
	\item There are isomorphisms 
	\[
	[\Rep^\square_{(E_\Db,\cE_\Db),\Ad} / Z(\cE_\Db)] \lrisom \Rep_{E_\Db,{D^u_{E_\Db}}} \lrisom \Rep_\Db
	\]
	of stacks on $\tfg_{\bZ_p}$. 
	
	\item Let $B$ be a commutative Noetherian local $\bZ_p$-algebra. Given an adapted representation $(E_\Db,\cE_\Db) \to M_d(B)$, the map $E_\Db \to M_d(B)$ is a compatible representation of $(E_\Db,D^u_{E_\Db})$.
	\item Let $(E, \rho, D)$ be a Cayley-Hamilton representation of $G$ with residual pseudorepresentation $\Db$. Then there is an $A$-GMA structure $\cE$ on $E_A$ such that $D_{\cE}=D$ and such that the map $(E_\Db,D_{\cE_\Db}) \to (E,D_{\cE})$ is adapted.
	\end{enumerate}
\end{cor}

\begin{rem}
\label{rem:formal adapted rep space}
Let $(E,\cE)$ be an $R$-GMA where $R \in \tfg_{\Z_p}$ is local. Let $S$ be the $R$-algebra from Theorem \ref{thm:Ad_main} so that $\Rep^\square_{(E,\cE),\Ad}=\Spec(S)$. Restricting the functor $\Rep^\square_{(E,\cE),\Ad}$  to the subcategory $\tfg_{R}$ of the category of $R$-algebras, we obtain an affine formal scheme $\Spf(\hat S)$, where $\hat S$ is the $\m_R S$-adic completion of $S$. We also denote $\Spf(\hat S)$ by $\Rep^\square_{(E,\cE),\Ad}$, abusing notation. This is how we consider $[\Rep^\square_{(E_\Db,\cE_\Db),\Ad} / Z(\cE_\Db)]$ as a formal stack on $\tfg_{\bZ_p}$.
\end{rem}

\begin{proof}
Statement (1)\ follows from Proposition \ref{prop:adapted_GMA} and Theorem \ref{thm:ED-compat}, while (2)\ follows from the statement of Theorem \ref{thm:CH-GMA} that $D^u_{E_\Db} = D_{\cE_\Db}$. Statement (3)\ follows from the second part of Theorem \ref{thm:CH-GMA}.
\end{proof}

\begin{lem}
\label{lem:Einj}
Assume that $\Db$ is multiplicity-free. Let $(E_\Db,D^u_{E_\Db}) \to (E,D)$  be a morphism of Cayley-Hamilton algebras, where $E$ is an $R$-algebra.
\begin{enumerate}
\item For any non-zero $x \in E$, there is a commutative local $R$-algebra $B$ of finite cardinality and a compatible representation $\rho_B: E \ra M_d(B)$ such that $\rho_B(x) \neq 0$.
\item Let $(E_\Db,D^u_{E_\Db}) \to (E',D')$ be another morphism of Cayley-Hamilton algebras, and assume that $E_\Db \to E$ and $E_\Db \to E'$ are surjective. If, for all $B$ as in (1) and all compatible representations $\rho_B:E_\Db \to M_d(B)$, the map $\rho_B$ factors through $E$ if and only if it factors through $E'$, then there is a canonical isomorphism $(E, D) \isoto (E', D')$ of Cayley-Hamilton algebras.
\end{enumerate}
\end{lem}
\begin{proof}
Since $\Db$ is multiplicity-free, we may fix a GMA structure $\cE_\Db$ on $E_\Db$ as in Theorem \ref{thm:CH-GMA}. By Corollary \ref{cor:compat=adapt}(3), this gives a GMA structure $\cE$ on $E$ such that $(E_\Db,D^u_{E_\Db}) \to (E,D)$ induces a morphisms of GMAs $(E_\Db, \cE_\Db) \to (E,\cE)$, and similarly for $(E',D')$.  By Corollary \ref{cor:compat=adapt}(3) we may work with adapted representations of these GMAs in the place of compatible representations of the Cayley-Hamilton algebras. We will do this for the remainder of the proof.

(1)\ By Theorem \ref{thm:Ad_main} and Remark \ref{rem:formal adapted rep space}, there is $\hat S \in \tfg_{R}$ such that $\Rep^\square_{(E,\cE),\Ad}=\Spf(\hat S)$ and, moreover, $E \rinj M_d(S)$ splits as an $R$-module map. Therefore $\rho_\Ad: E \to M_d(\hat S)$ remains injective.

Let $x \in E$ be a non-zero element, so $\rho_\Ad(x) \ne 0$. Let $y \in \hat S$ be a non-zero entry in the matrix $\rho_\Ad(x)$. By Lemma \ref{lem:tfg algebra}, there is a commutative local $R$-algebra $B$ of finite cardinality and an $R$-algebra homomorphism $f:\hat S \to B$ such that $f(y) \ne 0$. Then the composite $f \circ \rho_\Ad: E \to M_d(B)$ is an adapted representation such that $(f \circ \rho_\Ad)(x) \neq 0$.

(2)\ Since the maps $(E_\Db,D^u_{E_\Db}) \to (E,D)$ and $(E_\Db,D^u_{E_\Db}) \to (E',D')$ are morphisms of pseudorepresentations and the maps $g:E_\Db \to E$ and $g':E_\Db \to E'$ are surjective, it suffices to show $\ker(g)=\ker(g')$. Assume for a contradiction that there exists $x \in \ker(g)$ with $x \not \in \ker(g')$. Since $g'(x) \ne 0$, part (1)\ implies that there is a commutative local $R$-algebra $B$ of finite cardinality and a compatible representation $\rho'_B: E' \ra M_d(B)$ such that $\rho'_B(g'(x)) \neq 0$. The adapted representation $\tilde{\rho}_B=\rho'_B \circ g': E_\Db \to M_d(B)$ factors through $E'$, so it must factor through $E$ by assumption. This implies that $\tilde{\rho}_B = \rho_B \circ g$ for some adapted representation $\rho_B$ of $E$. But then 
\[
0=\rho_B(g(x)) =\tilde{\rho}_B(x)= \rho'_B(g'(x)) \neq 0
\]
a contradiction.
\end{proof}

\subsection{Condition $\cC$ in the residually multiplicity-free case}
\label{subsec:RMF-P}

Let $G$ be a profinite group satisfying the $\Phi_p$ finiteness condition. Fix a stable condition $\cC \subset \Mt$ as in Definition \ref{defn:stable} and a residual pseudorepresentation $\Db: G \to \bF$. By Theorem \ref{thm:univ_CH_P}, there is a universal Cayley-Hamilton algebra with condition $\cC$, denoted $(E^\cC_\Db, D_{E^\cC_\Db})$. In the case that $\Db$ is multiplicity-free, the following theorem gives an alternate characterization of $(E^\cC_\Db, D_{E^\cC_\Db})$.

\begin{thm}
\label{thm:GMA-P}
Assume that $\Db$ is multiplicity-free.
\begin{enumerate}
\item Let $B$ be a commutative local $R_\Db$-algebra of finite cardinality. Let $\rho_B: (E_\Db,D^u_{E_\Db}) \to M_d(B)$ be a compatible representation, and let $V_B \cong B^d$ denote the corresponding object of $\Mt$. Then $V_B$ satisfies $\cC$ if and only if $\rho_B$ factors through $E^\cC_\Db$. 
\item The property of (1) characterizes the quotient $E_\Db \to E^\cC_\Db$. 
\end{enumerate}
\end{thm}

\begin{proof}
Part (1) follows from Theorem \ref{thm:ModToCH} and Theorem \ref{thm:univ_CH_P}. Part (2) is immediate from Lemma \ref{lem:Einj}(2).
\end{proof}

When $\Db$ is not multiplicity-free, we only know how to characterize $E^\cC_\Db$ by Theorem \ref{thm:univ_CH_P}, cf.\ \cite[\S1.3.4]{BC2009}. 

This theorem furnishes a convenient way to prove that, in certain cases where a pseudorepresentation comes from a representation, property $\cC$ for the pseudorepresentation is related to property $\cC$ for a representation inducing it. 

\begin{cor}
\label{cor:tenauthors}
Let $\Db$ be multiplicity-free with $\bro_{\Db}^{ss} = \bro_1 \oplus \cdots \bro_r$ over $\F$, where the $\bro_i$ are absolutely irreducible and pairwise non-isomorphic of dimension $d_i$. For $i=1,\dots,r$, let $\rho_i:G \to M_{d_i}(A)$ be a deformation of $\bro_i$, and let $D=\psi(\rho_1 \oplus \cdots \oplus \rho_r)$. Then $D$ has $\cC$ if and only if there are deformations $\rho_i'$ of $\bro_i$ with $\cC$ for $i=1,\dots,r$ such that $D=\psi(\rho_1' \oplus \cdots \oplus \rho_r')$.
\end{cor}
\begin{proof}
If $D=\psi(\rho_1' \oplus \cdots \oplus \rho_r')$ for such $\rho'_i$, then its clear that $D$ has $\cC$. Now assume $D$ has $\cC$, so, by definition, there is a Cayley-Hamilton representation $(E,\rho,D')$ of $G$ that has $\cC$ such that $D= D'\circ \rho$. Since $\Db$ is multiplicity-free, we may assume that $D'=D_\cE$ for some GMA-structure $\cE=(\{e_i\},\{\phi_i\})$ on $E$ such that $\rho_i'=e_i\rho e_i: G \to M_{d_i}(A)$ is a deformation of $\bro_i$. Replacing $E$ by the image of $A[G]\to E$ if necessary, we may assume that the maps $\rho_{i,j}$ given by the GMA structure are surjective.

Now, the fact that $D=\psi(\rho_1 \oplus \cdots \oplus \rho_r)$ implies that $\phi_{i,j,k}=0$ for all triples $i,j,k$ of distinct integers with $1 \le i,j,k \le r$ (see \cite[\S 1.5.1]{BC2009}). Then the sum of projection maps
\[
E \to \bigoplus_{i=1}^r e_i E e_i =\bigoplus_{i=1}^r M_{d_i}(A).
\]
is an $A$-algebra homomorphism. The resulting map
\[
G \to E \to M_d(A)
\]
is $\rho':=\oplus_{i=1}^r \rho_i'$. Since $(E,D')$ has $\cC$, the map $E_\Db \to E$ factors through $E_\Db^\cC$. This implies that $\rho'$ has $\cC$ by Theorem \ref{thm:GMA-P}(1), and hence each $\rho_i'$ has $\cC$ by stability.
\end{proof}
\begin{rem}
The preceding corollary can be useful in proving automorphy lifting theorems. Indeed, it is a general version of an argument used in the proof of \cite[Prop.~4.4.6]{tenauthors}.
\end{rem}
\begin{warn}
When applying these results to the situation of \S\ref{subsec:global conditions} below (deformations of global Galois representations, where $\cC$ is a condition on the local representation), the preceding corollary only applies to deformations $D$ that are \textbf{globally} reducible. For example, take $\cC$ to be the condition ``de Rham at $p$" and $\Db=\psi(\F_p(1)\oplus \F_p):G_\Q \to \F_p$. Suppose $D=\psi(\rho)$ where $\rho:G_\Q \to \GL_2(\Z_p)$ is irreducible, but $\rho|_{G_p}$ is a non-split extension of $\Z_p(1)$ by $\Z_p$ (and hence not de Rham). Then $D|_{G_p}=\psi(\rho|_{G_p})=\psi(\Z_p(1) \oplus \Z_p)$ is a sum of two de Rham characters, and hence de Rham, but $D$ is not (globally) de Rham. 
\end{warn}

If $\bro$ is absolutely irreducible, let $R_\bro^\cC$ be the deformation-with-$\cC$ ring defined by Ramakrishna (denoted $R_\cC(\bro)$ in \cite[Prop.~1.2]{ramakrishna1993}) and let $\rho^\cC_\bro$ be the universal deformation.
\begin{cor}
\label{cor:us=Ram}
Suppose that $\bro_\Db^{ss}=\bro$ is absolutely irreducible. Then there is canonical isomorphism $R_\Db^\cC \to R_\bro^\cC$ determined by the pseudodeformation $\psi(\rho_\bro^\cC)$ of $\Db$.
\end{cor}
\begin{proof}
It is well-known that $R_\Db=R_\bro$ -- that is, every deformation $D$ of $\Db$ is of the form $D=\psi(\rho)$ for a unique deformation $\rho$ of $\bro$. (Using pseudo-characters, this is due to Carayol and Mazur. For pseudo-representations, this is proved by Chenevier \cite[Thm.\ 2.22(i)]{chen2014}; it also follows from Theorem \ref{thm:CH-GMA}.)  In this situation, by the previous corollary, $D$ has $\cC$ if and only if $\rho$ has $\cC$, so the two quotients $R_\Db^\cC$ and $R_\bro^\cC$ of $R_\Db=R_\bro$ coincide.
\end{proof}

This corollary shows that our $R_\Db^\cC$ is really a generalization of Ramakrishna's theory.

\section{Deformation conditions and extensions}
\label{sec:gma_ext}

In \cite{BC2009}, Bella\"{i}che and Chenevier used the explicit nature of the generalized matrix algebras discussed in the previous section to partially compute the universal pseudodeformation ring in terms of extension groups, in the category $\mathrm{Mod}_{\Z_p[G]}$, of the constituent irreducible representations of the residual representation. This calculation is a generalization of the description of the tangent space of the deformation ring of an irreducible representation in terms of adjoint cohomology.

In this section, we show that pseudodeformations-with-$\cC$ may be controlled by similar calculations, but where the extension groups are taken in the category $\cC$. Whereas Bella\"{i}che and Chenevier work with arbitrary multiplicity-free $\Db$, we, for reasons of simplicity and clarity, have chosen only to consider the case where $\Db$ is a sum of two characters (but see Remark \ref{rem:ExtP_expectation}).

In this situation, where $\Db = \psi(\bar\chi_1 \oplus \bar\chi_2)$ for distinct characters $\bar\chi_1$ and $\bar\chi_2$, we now outline what we can compute about $R_\Db^\cC$. First, we consider a simpler deformation ring $R_\Db^{\red, \cC}$ where we only consider deformations $D$ that remain reducible (i.e.~are of the form $D=\psi(\chi_1 \oplus \chi_2)$). This ring can be computed in terms of deformation rings of the characters $\bar\chi_1$ and $\bar\chi_2$. Next, we wish to compute the kernel $J^\cC=\ker(R_\Db^\cC \to R_\Db^{\red, \cC})$, which is called the reducibility ideal. By the work of Bella\"{i}che and Chenevier, this ideal can be described as a quotient of $B^\cC \otimes_{R_\Db^\cC} C^\cC$, where these modules come from the GMA structure on $E_\Db^\cC$. Finally, $B^\cC/J^\cC B^\cC$ and $C^\cC/J^\cC C^\cC$ can be described in terms of extension groups in the category $\cC$ of the universal deformations $\chi_1$ and $\chi_2$. In total, this gives an upper bound for the size of $J^\cC/(J^\cC)^2$ in terms of  extension groups of $\chi_1$ and $\chi_2$ in the category $\cC$, where $R_\Db^\cC/J^\cC$ is a quotient of $R_\Db^\cC$ that is fairly simple. 

\begin{rem}
\label{rem:use for R=T}
When proving modularity lifting theorems, this is often the kind of upper bound one needs. In that situation, one has a surjection $R_\Db^\cC \to \bT$, where $\bT$ is the quotient of `modular deformations.' To prove modularity, one needs to show that $R_\Db^\cC$ is not `too big' -- i.e.\ give an upper bound on the size of $R_\Db^\cC$. In this case, the extension groups in $\cC$ are given as certain Selmer groups in Galois cohomology (see Example \ref{eg:selmer H1}). Often, the size of these Selmer groups can be bounded in terms of modular invariants (e.g.~zeta values), and the resulting upper bound on $J^\cC/(J^\cC)^2$ can be enough to show that $R_\Db^\cC \to \bT$ is injective. This is the strategy of \cite{WWE3,WWE5}.
\end{rem}

\begin{rem}
In fact, our results allow a similar computation about $I/I^2$ for any ideal $I$ containing $J^\cC$. In particular, taking $I$ to be the maximal ideal of $R_\Db^\cC$, we get a description of the (co)tangent space of $R_\Db^\cC$. This computation about $J^\cC/(J^\cC)^2$ is giving information about the `relative tangent space' over the reducible deformation ring, and is, in practice, often more useful.
\end{rem}

We now outline the contents of this section. In \S \ref{subsec:RedExt}, we review the results of Bella\"{i}che--Chenevier. In \S \ref{subsec:PExt}, we prove our generalization to deformations-with-$\cC$.

\subsection{Conventions}
\label{subsec:ext-conventions}
 Throughout this section, we fix $G$, a profinite group satisfying the $\Phi_p$ finiteness condition, and two distinct characters $\bar\chi_1, \bar\chi_2: G \ra \bF^\times$. We let $\Db = \psi(\bar\chi_1 \oplus \bar\chi_2)$. 

By Theorem \ref{thm:CH-GMA}, we can and do fix a $R_\Db$-GMA structure $\cE_\Db =(\{e_1,e_2\},\{\phi_1,\phi_2\})$ on $E_\Db$. We write $(E_\Db,\cE_\Db)$ and $\rho^u: G \to E_\Db^\times$ as
\begin{equation}
\label{eq:CH-GMA}
E_\Db \cong \begin{pmatrix} R_\Db & B^u \\ C^u & R_\Db\end{pmatrix}, \qquad \rho^u(\sigma) = \ttmat{\rho^u_{1,1}(\sigma)}{\rho^u_{1,2}(\sigma)}{\rho^u_{2,1}(\sigma)}{\rho^u_{2,2}(\sigma)}
\end{equation}
and write $m: B^u \otimes_{R_\Db} C^u \to R_\Db$ for the induced map, as in Example \ref{eg:GMA1,1}. For $b \in B^u$ and $c \in C^u$, we define $b \cdot c = m(b \otimes c) \in R_\Db$. We can and do assume that the idempotents are ordered so that the image of $\rho_{i,i}(\sigma)$ under $R_\Db \to \bF$ is $\bar\chi_i(\sigma)$. 

By Corollary \ref{cor:compat=adapt}(3), a Cayley-Hamilton representation $(E, \rho, D)$ of $G$ with residual pseudorepresentation $\Db$  inherits a GMA structure from the data above. We will use matrix notation for $E$ and $\rho$ according to this structure, as in \eqref{eq:CH-GMA} (for example, writing the coordinates of $\rho$ as $\rho_{i,j}$).

\subsection{Review of reducibility, extensions, and GMAs}
\label{subsec:RedExt}

In this subsection, we review known results of \cite[\S1.5]{BC2009} relating the structure of GMAs to $\Ext$-groups. 

\begin{defn}
\label{defn:reducible}
Let $A \in \hat\cC_{W(\F)} $. We call a pseudodeformation $D: G \ra A$ of $\Db$ \emph{reducible} if $D=\psi(\chi_1 \oplus \chi_2)$ for characters $\chi_i: G \ra A^\times$ deforming $\bar\chi_i$. Otherwise, we call $D$ \emph{irreducible}. 

A GMA representation $(E, \rho: G \ra E^\times, D_\cE)$ of $G$ with residual pseudorepresentation $\Db$ is called \emph{reducible} (resp.\ \emph{irreducible}) provided that pseudodeformation $D_\cE \circ \rho : G \ra A$ is reducible (resp.\ irreducible). 
\end{defn}

\begin{prop}
\label{prop:reducible}
Let $A$ be a commutative Noetherian local $\bZ_p$-algebra. Let $D: G \ra A$ pseudorepresentation with residual pseudorepresentation $\Db$. 
\begin{enumerate}
\item $D$ is reducible if and only if $D=\psi(\rho)$ for some GMA representation $\rho$ with scalar ring $A$ such that $\rho_{1,2}(G) \cdot \rho_{2,1}(G)$ is zero.
\item Let $J = B^u \cdot C^u \subset R_\Db$ be the image ideal of $B^u \otimes_{R_\Db} C^u$ in $R_\Db$ under $m : B^u \otimes_{R_\Db} C^u \ra R_\Db$. Then $D$ is reducible if and only if the corresponding local homomorphism $R_\Db \ra A$ kills $J$. 
\end{enumerate}
\end{prop}
\begin{proof}
See \cite[\S1.5.1]{BC2009}. 
\end{proof}

In light of Proposition \ref{prop:reducible}, we establish the following terminology.
\begin{defn}
\label{defn:univ_red}
The ideal $J \subset R_\Db$ of Proposition \ref{prop:reducible}(2) is called the \emph{reducibility ideal} of $R_\Db$. The image of $J$ under the map $R_\Db \ra A$ corresponding to a pseudodeformation $D : G \ra A$ of $\Db$ is called the \emph{reducibility ideal of $D$}. 

We define $E^\red_\Db$ to be the Cayley-Hamilton quotient of $E_\Db$ by $J E_\Db$, which as in Example \ref{eg:CH-quot-with-scalar-ideal} is the usual algebra quotient $E^\red_\Db = E_\Db/JE_\Db$. Its scalar ring is $R^\red_\Db = R_\Db/J$ and is called the \emph{universal reducible pseudodeformation ring} for $\Db$. We let $(E^\red_\Db, \rho^\red, \cE^\red_\Db)$ denote the corresponding GMA representation of $G$, and write the decomposition arising from \eqref{eq:CH-GMA} as
\begin{equation}
\label{eq:red-CMA}
E^\red_\Db = \begin{pmatrix} R^\red_\Db & B^\red \\ C^\red & R^\red_\Db\end{pmatrix}.
\end{equation}
\end{defn}
\begin{warn}
Being reducible implies that the map $m: B^\red \otimes C^\red \to R^\red_\Db$ is zero, but it \textbf{does not} imply that either $B^\red$ or $C^\red$ is zero. In fact, for a fixed $g \in G$, it is possible that both $\rho^\red_{1,2}(g) \in B^\red$ and $\rho^\red_{2,1}(g) \in C^\red$ are non-zero.
\end{warn}

Let for $i=1,2$, let $R_{\bar\chi_i}$ denote the universal deformation ring of $\bar\chi_i$ and let $\chi_i^u: G \to R_{\bar\chi_i}^\times$ denote the universal deformation.

\begin{prop}
\label{prop:univ_red}
\hfill 
\begin{enumerate}
\item If $\rho:G \to (E,D)$ is a GMA representation of $G$ with scalar ring $A$ and residual pseudorepresentation $\Db$, then the resulting GMA map $(E_\Db, \cE_\Db) \ra (E,\cE)$ factors through $(E^\red_\Db, \cE^\red_{\Db})$ if and only if $\rho$ is reducible. 
\item There is a canonical isomorphism $R^\red_\Db \cong R_{\bar\chi_1} \hat\otimes_{W(\bF)} R_{\bar\chi_2}$. Letting $\chi_i = \chi_i^u \otimes_{R_{\bar\chi_i}} R^\red_\Db$, the universal reducible pseudodeformation of $\Db$ is $\psi(\chi_1\oplus \chi_2)$.
\item In terms of the decomposition \eqref{eq:red-CMA}, we have $\rho^\red_{i,i} = \chi_i$.
\end{enumerate}
\end{prop}

\begin{proof}
(1)\ By Lemma \ref{lem:max_CH_quot}, the GMA map $(E_\Db, \cE_\Db) \ra (E,\cE)$ factors through $(E^\red_\Db, \cE^\red_{\Db})$ if and only if $E_\Db \to E$ sends $JE_\Db$ to zero. Because $R \ra E$ is injective whenever there is a pseudorepresentation $D: E \ra R$ \cite[Lem.\ 5.2.5]{WWE1}, $JE_\Db$ maps to zero in $E$ if and only if $R_\Db \to A$ factors through $R_\Db^\red$, which, by Proposition \ref{prop:reducible}(2), is equivalent to $D \circ \rho$ being reducible. 

(2)\ By Yoneda's lemma, it suffices to construct a canonical functorial isomorphism 
\[
\Hom(R^\red_\Db,A) \cong \Hom(R_{\bar\chi_1} \hat\otimes_{W(\bF)} R_{\bar\chi_2},A)
\]
for $A \in \hat\cC_{W(\bF)}$. Given $R_{\bar\chi_1} \hat\otimes_{W(\bF)} R_{\bar\chi_2} \to A$, we define a reducible deformation $D$ of $\Db$ by $D=\psi(\chi_1^u \otimes A \oplus \chi_2^u \otimes A)$, which determines an element of $\Hom(R^\red_\Db,A)$ by Proposition \ref{prop:reducible}(2) and the universal property of $R_\Db$. 

Conversely, given $R^\red_\Db \to A$, consider the GMA representation $\rho^\red \otimes_{R^\red_\Db} A$. Since $B^\red \cdot C^\red=0$, we have that $\rho^\red_{i,i} \otimes_{R^\red_\Db} A: G \to A^\times$  (for $i=1,2$) is a character that, by the conventions of \S \ref{subsec:ext-conventions}, is a deformation of $\bar\chi_i$. This pair of characters determines an element of $\Hom(R_{\bar\chi_1} \hat\otimes_{W(\bF)} R_{\bar\chi_2},A)$. 

We have defined maps between $\Hom(R^\red_\Db,A)$ and $\Hom(R_{\bar\chi_1} \hat\otimes_{W(\bF)} R_{\bar\chi_2},A)$. The reader can check these are mutually inverse and functorial in $A$. 

(3)\ It follows from the construction of the isomorphism in (2).
\end{proof}

The following key result relates $R^\red_\Db[G]$-module extensions to the structure of $E^\red_\Db$. For ease of notation, we write $\chi_i$ for the base change from $R_{\bar\chi_i}$ to $R^\red_\Db$ of the universal deformation $\chi_i^u$ of $\bar\chi_i$. 

\begin{prop}
\label{prop:GMA-Ext}
Let $A \in \hat\cC_{W(\F)}$ and let $M$ be a finitely generated $A$-module. For $i=1,2$, let $\chi_{i,A}:G \to A^\times$ be characters deforming $\bar\chi_i$. By Proposition \ref{prop:univ_red}, this induces a unique local homomorphism $R_\Db^\red \to A$. There is a natural isomorphism  
\[
\Hom_{A}(B^\red\otimes_{R^\red_\Db} A,M) \lrisom \Ext^1_{A[G]}(\chi_{2,A}, \chi_{1,A} \otimes_{A} M). 
\]
as well as a similar isomorphism in the $C$-coordinate. 
\end{prop}

\begin{proof}
The details may be found in \cite[Lem.\ 4.1.5, proof of (4.1.7)]{WWE2}, cf.~also \cite[Thm.\ 1.5.6]{BC2009}. We reproduce the construction of the map here with notation that will be convenient in \S\ref{subsec:PExt}.   

Let $E_M=\sm{A}{M}{0}{A}$, with GMA structure as in Example \ref{eg:m0}. Given a homomorphism $f:B^\red \otimes_{R^\red_\Db} A \ra M$, we have morphism of GMAs $E^\red_\Db\otimes_{R^\red_\Db} A \to E_M$ as the composition
\begin{equation}
\label{eq:EM-factor}
E^\red_\Db\otimes_{R^\red_\Db} A \xrightarrow{\sm{\mathrm{id}}{\mathrm{id}}{0}{\mathrm{id}}} \ttmat{A}{B^\red \otimes_{R^\red_\Db} A}{0}{A} \xrightarrow{\sm{\mathrm{id}}{f}{0}{\mathrm{id}}} \ttmat{A}{M}{0}{A}.
\end{equation}
Using the fact that $B^\red \cdot C^\red =0$, we see that $e_1E_Me_2 =  \sm{0}{M}{0}{0}$ is a left $E_\Db$-submodule of $E_M e_2 = \sm{0}{M}{0}{A}$. Noting that $e_1E_Me_2 \cong \chi_{1,A} \otimes_{A} M$ and $E_Me_2 / e_1E_Me_2 \cong \chi_{2,A}$ as $A[G]$-modules, we obtain an exact sequence
\begin{equation}
\label{eq:GMA-SES}
0 \lra \chi_{1,A} \otimes_{A} M \ra E_M e_2 \lra \chi_{2,A} \lra 0,
\end{equation}
which determines the corresponding element of $\Ext^1_{A[G]}(\chi_{2,A}, \chi_{1,A} \otimes_{A} M)$. Conversely, any such extension can be realized in the form $E_Me_2$. 
\end{proof}

\subsection{GMA structures corresponding to extensions with an abstract property}
\label{subsec:PExt}
Let $\cC \subset \Mt$ be a stable property, as in Definition \ref{defn:stable}. The following lemma is a well-known consequence of stability, and we leave the proof to the reader. To state the lemma clearly, we introduce some notation for extension classes. If $E$ is an algebra and $V_1, V_2$ are $E$-modules, given an extension class $c \in \Ext^1_{E}(V_1,V_2)$, and an exact sequence
\begin{equation}
\label{eq:ses-representing-ext-class}
0 \to V_2 \to V \to V_1 \to 0
\end{equation}
representing $c$, we call $V$ an \emph{extension module} for $c$.

Let $V_1, V_2 \in \cC$, and assume that $V_1, V_2 \in \mathrm{Mod}_{R[G]}$ as well for commutative $\Z_p$-algebra $R$. Define $\Ext^1_{R[G],\cC}(V_1,V_2)$ as the subset of $\Ext^1_{R[G]}(V_1,V_2)$ consisting of extension classes $c$ such that, for every extension module $V$ for $c$, we have $V \in \cC$. 

\begin{lem}
\label{lem:C-exts}
With $V_1, V_2$ as above, we have the following.
\begin{enumerate}
\item For a class $c \in \Ext^1_{R[G]}(V_1,V_2)$, if some extension module $V$ for $c$ has $V \in \cC$, then $c \in \Ext^1_{R[G],\cC}(V_1,V_2)$. 
\item The subset $\Ext^1_{R[G],\cC}(V_1,V_2) \subset \Ext^1_{R[G]}(V_1,V_2)$ is a sub-$R$-module.
\end{enumerate}
\end{lem}

\begin{eg}
\label{eg:selmer H1}
 Let $H_1,\dots,H_n \subset G$ be subgroups, and, for $i=1,\dots,n$, let $\cC_i \subset \mathrm{Mod}_{\Z_p[H_i]}^{\mathrm{fin}}$ be a stable condition. Assume that the condition $\cC \subset \Mt$ arises from the $\cC_i$ as in Example \ref{eg:fiber product of stable is stable}. Then $\Ext^1_{R[G],\cC}(V_1,V_2)$ is the kernel of the map
\[
\Ext^1_{R[G]}(V_1,V_2)  \to \bigoplus_{i=1}^n \frac{\Ext^1_{R[H_i]}(V_1,V_2)}{\Ext^1_{R[H_i],\cC_i}(V_1,V_2)}
\]
given by the restrictions $\Ext^1_{R[G]}(V_1,V_2) \to \Ext^1_{R[H_i]}(V_1,V_2)$ followed by the quotients. (This is sometimes referred to as a Selmer group.)
\end{eg}

Let $E^\cC_\Db$ be as in Theorem \ref{thm:GMA-P}, and let $E^{\cC,\red}_\Db=E^\cC_\Db / J E^\cC_\Db$ where $J \subset R_\Db$ is the reducibility ideal. Following the notation of \eqref{eq:CH-GMA}, we write them as 
\[
E^\cC_\Db =\begin{pmatrix} R^\cC_\Db & B^\cC \\ C^\cC & R^\cC_\Db\end{pmatrix}, \qquad E^{\cC,\red}_\Db = \begin{pmatrix} R^{\cC,\red}_\Db & B^{\cC,\red} \\ C^{\cC,\red} & R^{\cC,\red}_\Db\end{pmatrix}.
\]
We denote the Cayley-Hamilton representations by $\rho^\cC:G \to (E^\cC_\Db)^\times$ and $\rho^{\cC,\red}:G \to (E^{\cC,\red}_\Db)^\times$.

By Ramakrishna's result \cite[Prop.\ 1.2]{ramakrishna1993}, for $i=1,2$, there is a quotient $R_{\bar\chi_i} \onto R_{\bar\chi_i}^\cC$ that represents the functor of deformations of $\bar\chi_i$ having property $\cC$ and $\chi_i^{u,\cC} = \chi_i^{u} \otimes_{R_{\bar\chi_i}} R_{\bar\chi_i}^\cC$ is the universal deformation with property $\cC$. 

\begin{prop}
\label{prop:univ-red-with-C}
There is a canonical commutative diagram 
\[
\xymatrix{
R^{\cC,\red}_\Db \ar[r]^-\sim &  R^\cC_{\bar\chi_1} \hat\otimes_{W(\bF)} R^\cC_{\bar\chi_2} \\
R^\red_\Db \ar@{->>}[u] \ar[r]^-\sim & R_{\bar\chi_1} \otimes_{W(\bF)} R_{\bar\chi_2} \ar@{->>}[u]
}
\]
where the vertical maps are the quotients and the lower horizontal map is the isomorphism of Proposition \ref{prop:univ_red}(2).
\end{prop}

\begin{proof}
For simplicity of notation, let $R=R^\cC_{\bar\chi_1} \hat\otimes_{W(\bF)} R^\cC_{\bar\chi_2}$ and, for $i=1,2$, let $\chi_i^\cC = \chi_i^{u,\cC} \otimes_{R^\cC_{\bar\chi_i}} R$. 

By Proposition \ref{prop:univ_red}(2), the composite map $R^\red_\Db \to R$ corresponds to the reducible pseudodeformation $D=\psi(\chi_1^\cC \oplus \chi_2^\cC)$. The $G$-module $N=\chi_1^\cC \oplus \chi_2^\cC$ is a faithful Cayley-Hamilton $G$-module with Cayley-Hamilton algebra $(E=\sm{R}{0}{0}{R}, D_\cE)$. By Theorem \ref{thm:ModToCH}, the map $(E_\Db,D^u_{E_\Db}) \to (E,D_\cE)$ factors through $E_\Db^\cC$. Since $D$ is reducible, Proposition \ref{prop:reducible} implies that it further factors through $E_\Db^{\cC,\red}$. This implies that $R^\red_\Db \to R$ factors through $R^{\cC,\red}_\Db$.

On the other hand, the composite map $R_{\bar\chi_1} \otimes_{W(\bF)} R_{\bar\chi_2} \to R^{\cC,\red}_\Db$ corresponds to the pair of characters $e_1\rho^{\red,\cC}e_1, e_2\rho^{\red,\cC}e_2: G \to  R^{\cC,\red}_\Db$. Since these characters are quotient $R^\red_\Db[G]$-modules of the $R^\red_\Db[G]$-module $E^{\cC,\red}_\Db$, which has $\cC$, we see that they both have $\cC$ as well. This implies that the map $R_{\bar\chi_1} \otimes_{W(\bF)} R_{\bar\chi_2} \to R^{\cC,\red}_\Db$ factors through $R$, completing the proof.
\end{proof}

We will write $\chi^\cC_i$ for the base change from $R^\cC_{\bar\chi_i}$ to $R^{\cC,\red}_\Db$ of the universal deformation of $\bar\chi_i$ satisfying $\cC$. 

\begin{thm}
	\label{thm:BCs and exts}
Let $A \in \hat\cC_{W(\F)}$ and let $M$ be a finitely generated $A$-module. For $i=1,2$, let $\chi_{i,A}:G \to A^\times$ be characters deforming $\bar\chi_i$ and having property $\cC$. By Proposition \ref{prop:univ-red-with-C}, this induces a unique local homomorphism $R_\Db^{\red,\cC} \to A$. There is a natural isomorphism  
\[
\Hom_{A}(B^{\cC,\red}\otimes_{R^{\red,\cC}_\Db} A,M) \lrisom \Ext^1_{A[G],\cC}(\chi_{2,A}, \chi_{1,A} \otimes_{A} M), 
\]
as well as a similar isomorphism in the $C$-coordinate. 	
\end{thm}

\begin{rem}
\label{rem:ExtP_expectation}
Throughout this section, we have restricted our attention to 2-dimensional pseudorepresentations and GMAs of type (1,1). However, in \cite[Thm.\ 1.5.6]{BC2009}, Bella\"{i}che and Chenevier prove a version of Proposition \ref{prop:GMA-Ext} for GMAs of any dimension $d$ and any type $(d_1,\dots,d_r)$. 

We have made this choice strictly for clarity of exposition. We fully expect a version of Theorem \ref{thm:BCs and exts} to be true for GMAs of any dimension $d$ and any type $(d_1,\dots,d_r)$, and that a proof can be given by applying the methods of this section to the general setup of \cite[Thm.\ 1.5.6]{BC2009}. For example, when $A = \F$, such a statement can be deduced directly from \cite[Thms.\ 11.2.1 and 12.3.1]{CarlAInf}. 

In fact, a version for types $(d_1,d_2)$ can be deduced from Theorem \ref{thm:BCs and exts} by Morita equivalence, as in \cite[\S1.3.2]{BC2009}. 
\end{rem}

\begin{proof}
We set $B=B^{\red}\otimes_{R^{\red}_{\Db}}A$ and $B^\cC = B^{\cC,\red}\otimes_{R^{\cC,\red}_\Db} A$ to simplify notation.  We consider the diagram
\[\xymatrix{
\Hom_{A}(B^{\cC},M) \ar@{^(->}[d] \ar@{-->}[r]^-\sim & \Ext^1_{A[G],\cC}(\chi_{2,A}, \chi_{1,A} \otimes_{A} M) \ar@{^(->}[d] \\
\Hom_{A}(B ,M) \ar[r]_-\sim^-\Psi & \Ext^1_{A[G]}(\chi_{2,A}, \chi_{1,A} \otimes_{A} M),
}\]
where $\Psi$ isomorphism of Proposition \ref{prop:GMA-Ext}, and where the dotted arrow is the isomorphism we wish to construct. Since this diagram is canonically isomorphic to the one obtained by replacing $A$ by the image of $A \to \End_A(M)$, we can and do assume that $M$ is a faithful $A$-module. Given a class $c \in \Ext^1_{A[G]}(\chi_{2,A}, \chi_{1,A} \otimes_{A} M)$, we have to show that the map $f_c=\Psi^{-1}(c):B  \to M$ factors through $B^{\cC}$ if and only if $c \in \Ext^1_{A[G],\cC}(\chi_{2,A}, \chi_{1,A} \otimes_{A} M)$.

Following the proof of Proposition \ref{prop:GMA-Ext}, we see that the class $c$ has an extension module $E_Me_2$ where $E_M =\sm{A}{M}{0}{A}$ is a GMA and $f_c$ induces a morphism of GMAs $E^{\red}_{\Db}\otimes_{R^{\red}_{\Db}}A \to E_M$ by
\begin{equation}
\label{eq:ext-to-map-of-GMA}
E^{\red}_{\Db}\otimes_{R^{\red}_{\Db}}A \xrightarrow{\sm{\mathrm{id}}{\mathrm{id}}{0}{\mathrm{id}}} \ttmat{A}{B}{0}{A} \xrightarrow{\sm{\mathrm{id}}{f_c}{0}{\mathrm{id}}} \ttmat{A}{M}{0}{A}.
\end{equation}
Since $M$ is a faithful $A$-module, we see that $E_Me_2$ is a faithful Cayley-Hamilton $G$-module with Cayley-Hamilton algebra $(E_M,D_{\cE_M})$. By Theorem \ref{thm:ModToCH}, $E_Me_2$ has $\cC$ if and only if $(E_M,D_{\cE_M})$ has $\cC$ as a Cayley-Hamilton representation of $G$. By Lemma \ref{lem:C-exts}, Lemma \ref{lem:max_CH_quot}, and Theorem \ref{thm:univ_CH_P}, we are reduced to showing that $f_c$ factors through $B^{\cC}$ if and only if the map $E_\Db \onto E^{\red}_{\Db}\otimes_{R^{\red}_{\Db}}A\to E_M$ given by \eqref{eq:ext-to-map-of-GMA} factors through $E_\Db^\cC$.

If $f_c$ factors through $B^{\cC}$, then we see that the map $E^{\red}_{\Db}\otimes_{R^{\red}_{\Db}}A\to E_M$ of \eqref{eq:ext-to-map-of-GMA} agrees with the map
\[
E^{\red}_{\Db}\otimes_{R^{\red}_{\Db}}A  \onto E^{\cC,\red}_{\Db}\otimes_{R^{\cC,\red}_{\Db}}A \xrightarrow{\sm{\mathrm{id}}{\mathrm{id}}{0}{\mathrm{id}}} \ttmat{A}{B^\cC}{0}{A} \xrightarrow{\sm{\mathrm{id}}{f_c}{0}{\mathrm{id}}} \ttmat{A}{M}{0}{A}.
\]
Hence the map $E_\Db \to E_M$ factors though $E^{\cC,\red}_{\Db}$, which is a quotient of $E_\Db^\cC$.

Conversely, suppose that $(E_\Db,D^u_{E_\Db}) \to (E_M,D_{\cE_M})$ factors through $(E_\Db^\cC,D_{E^\cC_\Db})$. Since $(E_M,D_{\cE_M})$ is reducible, this implies that the map $E_\Db \to E_M$ factors through a map $g:E^{\cC,\red}_{\Db}\otimes_{R^{\cC,\red}_{\Db}}A  \to E_M$. By \eqref{eq:ext-to-map-of-GMA}, we see that $f_c$ factors through $e_1ge_2: B^\cC \to M$.
\end{proof}

\section{Examples} 
\label{sec:ex}

The follow examples of conditions $\cC$ could be useful in arithmetic applications. 

\subsection{Unramified local condition} 
Let $\ell$ be a prime number and let $K$ be a finite field extension of $\bQ_\ell$. Let $G = \Gal(\overline{K}/K)$, which satisfies the $\Phi_p$ finiteness condition because it is topologically finitely generated. Let $H = I_\ell \subset G$ be the inertia subgroup. The condition $\mathrm{Mod}_{\Z_p[G/H]}^\mathrm{fin} \subset \Mt$, as in Example \ref{eg:factors through is stable}, is called \emph{unramified}. 

\subsection{Local conditions at $\ell=p$}
\label{subsec:Lex}

Retain the same notation as the previous subsection, but now assume $\ell=p$. Let $O_K \subset K$ denote the ring of integers. In this case, there are many conditions on representations of $G$ in $\Q_p$-vector spaces, coming from $p$-adic Hodge theory. Some of these conditions descend to $\Mt$ as follows.

\begin{defns} Let $V$ be an object of $\Mt$ and let $a \le b$ be integers.
\begin{enumerate} 
\item We call $V$ \emph{finite-flat} if there is a finite flat group scheme $\cG$ over $O_K$ such that $V \cong \cG(\overline{K})$ as $\Z_p[G]$-modules.
\item We call an object $V \in \Mt$ \emph{torsion crystalline (resp.~semi-stable) with Hodge-Tate weights in $[a,b]$} if there is a crystalline (resp.~semi-stable) representation $\rho:G \to \GL(W)$ with Hodge-Tate weights in $[a,b]$ and a $G$-stable $\Z_p$-lattice $T \subset W$ such that $V$ is isomorphic to a subquotient of $T$. 
\end{enumerate}
\end{defns}

Let $\cC_{\fl} \subset \Mt$ denote the full subcategory of finite-flat objects. Ramakrishna has proven that $\cC_{\fl}$ is stable \cite[\S2]{ramakrishna1993}. 

Let $\cC_{\mathrm{crys},[a,b]},  \cC_{\mathrm{st},[a,b]} \subset \Mt$ denote the full subcategories that are torsion crystalline (resp.\ semi-stable) with Hodge-Tate weights in $[a,b]$. Since the category of crystalline (resp.\ semi-stable) representations with Hodge-Tate weights in $[a,b]$ is closed under finite direct sums, we see that $\cC_{\mathrm{crys},[a,b]}$ and $\cC_{\mathrm{st},[a,b]}$ are closed under finite direct sums. They are also closed under isomorphisms and subquotients by definition, so we see that they are stable. 

\begin{rem}
It is known that, for a $\Z_p[G]$-module $M$ that is finitely generated and free as a $\Z_p$-module, $M \otimes_{\Z_p} \Q_p$ is crystalline (resp.~semi-stable) with Hodge-Tate weights in $[a,b]$ if and only if $M \otimes_{\Z_p} \Z/p^n\Z$ is torsion crystalline (resp.~semi-stable) with Hodge-Tate weights in $[a,b]$ for all $n \ge 1$. This was conjectured by Fontaine, and proven by Liu using results of Kisin, following partial results of Ramakrishna, Berger, and Breuil. It is also known that there is an equivalence of categories $\cC_\fl \cong \cC_{\mathrm{crys}, [0,1]}$ for $p > 2$. See \cite{liu2007} and the references given there.
\end{rem}

\begin{rem}
If $\Db$ is residually multiplicity-free and $\cC$ is one of the conditions about, then a Zariski-closed substack $\widetilde{\Rep}^\cC_\Db \subset \Rep_\Db$ and a quotient $R_\Db \onto \tilde{R}^\cC_\Db$ were constructed in \cite[\S7.1]{WE2018}. The two constructions will have the same generic fiber (over $\Z_p$), as this is where the construction of crystalline loci is done in \textit{loc.\ cit.} That is, $\tilde{R}^\cC_\Db[1/p] \cong R^\cC_\Db[1/p]$, and similarly for $\widetilde{\Rep}^\cC_\Db$ and $\Rep^\cC_\Db$. In particular, the geometric properties of the generic fibers proved in Prop.\ 6.4.4 and Cor.\ 7.1.5 of \emph{loc.\ cit.}\ apply to the generic fibers of the $\Rep^\cC_\Db$ and $R^\cC_\Db$ constructed in this paper. 

When the $\Z_p$-integral structures differ, we believe that the constructions in this paper will be better behaved: here, we make natively integral constructions, while the integral model of \textit{loc.\ cit.}\ is made by flat closure relative to $\Z_p$. 
\end{rem}

\subsection{Global conditions}
\label{subsec:global conditions}
Let $F$ be a number field with algebraic closure $\overline{F}$ and let $S$ be a finite set of places of $F$. Let $G=\Gal(F_S/F)$, where $F_S \subset \overline{F}$ is the maximal extension of $F$ unramified outside $S$. Then $G$ satisfies $\Phi_p$ by class field theory. 

For each $v \in S$, choose a decomposition group $G_v \subset G$ (so $G_v$ is of the type considered in the previous subsections) and a stable condition $\cC_v \subset \mathrm{Mod}_{\Z_p[G_v]}^\mathrm{fin}$. Then there is a corresponding stable condition $\cC \subset \Mt$, as in Example \ref{eg:fiber product of stable is stable}. The Selmer groups of Example \ref{eg:selmer H1} correspond to such a $\cC$.

\bibliographystyle{alpha}
\bibliography{CWEbib-2018-PG}

\enddocument